\documentclass{amsart}
\usepackage{amssymb}

\newtheorem{thm}{Theorem}[section]

\newtheorem{lem}[thm]{Lemma}
\newtheorem{rema}[thm]{Remark}

\theoremstyle{remark}

\numberwithin{equation}{section}

\newcommand{\Z}{\mathbb{Z}}
\newcommand{\C}{\mathbb{C}}
\begin{document}

\title[Twisted modules for N=2 VOSAs]{Twisted modules for N=2 supersymmetric vertex operator superalgebras arising from finite automorphisms of the N=2 Neveu-Schwarz algebra}

\author{Katrina Barron}
\address{Department of Mathematics, University of Notre Dame,
Notre Dame, IN 46556}
\email{kbarron@nd.edu}
\thanks{Supported in part by NSA grant MSPF-07G-169, and by a research grant from the Max Planck Institute for Mathematics, Bonn, Germany.}

\date{April 14, 2013}


\begin{abstract}
Twisted modules for N=2 supersymmetric vertex operator superalgebras are classified for the vertex operator superalgebra automorphisms which are lifts of a finite automorphism of the N=2 Neveu-Schwarz Lie superalgebra representation.    These include the Ramond-twisted sectors and mirror-twisted sectors for N=2 vertex operator superalgebras, as well as twisted modules related to more general ``spectral flow" representations of the N=2 Neveu-Schwarz algebra.   We present the construct of  twisted modules for free N=2 vertex operator superalgebras for all of the N=2 Neveu-Schwarz Lie superalgebra automorphisms of finite order.   We show how to extend these to lattice N=2 vertex operator superalgebras.  As a consequence, we also construct the Ramond-twisted sectors for free and lattice N=1 supersymmetric vertex operator superalgebras.  We show that the lifting of the mirror automorphism for the N=2 Neveu-Schwarz algebra to an N=2 vertex operator superalgebra is not unique and that different mirror map vertex operator superalgebra automorphisms of an N=2 vertex operator superalgebra can lead to non-isomorphic categories of mirror-twisted modules, as in the case of free and lattice N=2 vertex operator superalgebras.  
\end{abstract}

\maketitle

\tableofcontents

\section{Introduction and preliminaries}

We study twisted modules for  N=2 superconformal vertex operator superalgebras for vertex operator superalgebra automorphisms that arise from Virasoro-preserving automorphisms of the underlying N=2 Neveu-Schwarz algebra.  In particular, we characterize all such twisted modules in terms of the resulting representations of N=2 superconformal algebras.  We present explicit examples for free and lattice N=2 vertex operator superalgebras and show how to extend these constructions to lattice N=2 vertex operator superalgebras.  

If $g$ is an automorphism of a vertex operator superalgebra, $V$, then we have the notion of ``$g$-twisted $V$-module".  Twisted vertex operators were discovered and used in \cite{LW}.  Twisted modules for vertex operator algebras arose in the work of  I. Frenkel, J. Lepowsky and A. Meurman \cite{FLM1}, \cite{FLM2}, \cite{FLM3} in the course of the construction of the moonshine module vertex operator algebra. This structure came to be understood as an ``orbifold model" in the sense of conformal field theory and string theory.  Twisted modules are the mathematical counterpart of ``twisted sectors", which are the basic building blocks of orbifold models in conformal field theory and string theory.  The notion of twisted module for vertex operator superalgebras was developed in \cite{Li-twisted}.  In general, given a vertex operator algebra (let alone a vertex operator superalgebra) $V$ and an automorphism $g$ of $V$, it is an open problem as to how to
construct a $g$-twisted $V$-module.

An automorphism $g$ of a vertex operator superalgebra (VOSA), in particular, fixes the Virasoro vector, and thus also fixes the corresponding endomorphisms giving the representation of the Virasoro algebra.  For example, any VOSA, $V$, with $\mathbb{Z}_2$-grading given by $V^{(0)} \oplus V^{(1)}$, we have the {\it parity automorphism} $\sigma : v \mapsto (-1)^{|v|} v$ where $|v| = j$ if $v \in V^{(j)}$.

A VOSA is said to be ``N=1 or N=2 supersymmetric", if in addition to being a positive energy representation for the Virasoro algebra, it is a representation of the N=1 or N=2 Neveu-Schwarz algebra, respectively; see, for instance, \cite{DPZ1985}, \cite{B-announce}, \cite{B-vosas}, \cite{S}, \cite{B-n2axiomatic}.    The group of automorphisms of the N=1 Neveu-Schwarz algebra which preserve the Virasoro algebra is $\mathbb{Z}_2$, and is generated by the parity automorphism, $\sigma$ on the $\mathbb{Z}_2$-grading of the Lie superalgebra structure of the N=1 Neveu-Schwarz algebra.  
The group of automorphisms of the N=2 Neveu-Schwarz algebra over $\mathbb{C}$, which preserve the Virasoro algebra, is isomorphic to $\mathbb{C}^\times \times \mathbb{Z}_2$.  It is generated by a continuous family of automorphisms, denoted by $\sigma_\xi$ for $\xi \in \mathbb{C}^\times$, and an order two automorphism $\kappa$ called the ``mirror map".   If $\xi$ is a root of unity, then $\sigma_\xi$ is of finite order.  If $\xi = -1$, then $\sigma_{-1}$ is the parity map $\sigma$.   

Given an N=2 supersymmetric VOSA, $V$, some questions naturally arise:  When does $\kappa$ or $\sigma_\xi$ for $\xi \neq -1$, lift to an automorphism of $V$, and when is this lift unique?  When such an automorphism of the N=2 Neveu-Schwarz algebra does lift to an automorphism $g$ of $V$, what is the structure of a $g$-twisted $V$ module?  In this paper, we fully answer the second question, and we answer the first question for free and lattice N=2 VOSAs.

If the mirror automorphism $\kappa$ of the N=2 Neveu-Schwarz algebra lifts to a VOSA automorphism of an N=2 VOSA, $V$, then a ``mirror-twisted $V$-module" is naturally a representation of what we call the ``mirror-twisted N=2 superconformal algebra", which is also referred to as the ``twisted N=2 superconformal algebra" \cite{SS}, \cite{DG2001}, \cite{LSZ2010}, or the ``topological N=2 superconformal algebra"  \cite{Gato-Rivera2002}.  If the automorphism $\sigma_\xi$ of the N=2 Neveu-Schwarz algebra, for $\xi$ a root of unity, lifts to an VOSA automorphism of $V$, then we show that a ``$\sigma_\xi$-twisted $V$-module" is naturally a representation of one of the algebras in the one-parameter family of Lie superalgebras we call ``shifted N=2 superconformal algebras".   If $\xi = -1$, then $\sigma_\xi$ is the parity automorphism on the N=2 Neveu-Schwarz algebra, and a lift of $\sigma_\xi = \sigma_{-1}$ to an automorphism of an N=2 VOSA, $V$, always exists, namely as $\sigma$, the parity automorphism on $V$ as a VOSA, although such a lift is not necessarily unique, as we shall show.  The ``shifted N=2 superconformal algebra" resulting in this case of twisting by a lift of $\sigma_\xi = \sigma_{-1}$, is also called the ``N=2 Ramond algebra".   The N=2 Ramond algebra and the other shifted N=2 algebras are isomorphic, as Lie superalgebras, to the N=2 Neveu-Schwarz algebra via the ``spectral flow" operators, as was first realized in \cite{SS}.  The mirror-twisted N=2 algebra is not isomorphic to the N=2 Neveu-Schwarz algebra.   

The representation theory of the N=2 Neveu-Schwarz algebra has been studied in, for instance, \cite{DPZ1985}, \cite{DPY1986}, \cite{BFK1986}, \cite{Nam1986}, \cite{Dobrev1987}, \cite{Matsuo1987}, \cite{Kiritsis1988}, \cite{Dorrzapf1995}, \cite{ST1998}, \cite{FST1998}, \cite{FSST1999}, \cite{DG1999}, \cite{STF} and from a VOSA theoretic point of view in  \cite{Adamovic1999}, \cite{Adamovic2001}.  The representation theory of the N=2 Ramond algebra has been studied in, e.g, \cite{Dobrev1987},  \cite{ST1998}, \cite{FST1998}, \cite{Gato-Rivera2002}, \cite{FJS2007}, and of the mirror-twisted N=2 superconformal algebra in, e.g,  \cite{Dobrev1987}, \cite{DG2001}, \cite{Gato-Rivera2002}, \cite{IK2009}, \cite{LSZ2010}.

The realization of the N=2 Ramond algebra and the mirror-twisted N=2 superconformal algebra as arising from twisting an N=2 VOSA (or comparable structure) has long been known, e.g. \cite{SS}, \cite{BFK1986}, \cite{DPZ1985}.  However to our knowledge, the other algebras related to the N=2 Neveu-Schwarz algebra---the shifted N=2 superconformal algebras other than the N=2 Ramond algebra---have only been studied through the spectral flow operators (which do not preserve the Virasoro algebra).  We believe that the realization of these algebras as arising naturally as twisted modules for an N=2 VOSA is new. 

Thus this complete classification of the form of the twisted modules for an N=2 VOSA (in terms of representations of the various N=2 superconformal algebras) for finite automorphisms arising from Virasoro-preserving automorphisms of the N=2 Neveu-Schwarz algebra provides a uniform way of understanding and studying all of the N=2 superconformal algebras---the continuous one-parameter family of shifted N=2 Neveu-Schwarz algebras and the mirror-twisted N=2 superconformal algebra---in the context of the theory of VOSAs and their twisted modules.  

For all of these types of twisted modules arising from finite automorphisms of the N=2 Neveu-Schwarz algebra, we construct examples.  For the automorphisms $\sigma_\xi$, we indicate how these extend to free and lattice N=2 VOSAs and construct the corresponding twisted N=2 VOSA modules.  We also give the graded dimensions for these constructions.   This includes constructing and classifying the Ramond twisted sectors for both N=1 supersymmetric VOSAs and N=2 supersymmetric VOSAs in the case of free N=1 (resp. N=2) supersymmetric VOSAs and showing how these can be extended to lattice N=2 VOSAs.

For free and lattice N=2 VOSAs, we find two distinct mirror maps which represent distinct conjugacy classes within the group of automorphisms for the given N=2 VOSA and thus result in non-isomorphic categories of mirror-twisted modules.  We carry out the construction of the mirror-twisted modules for one of these mirror maps and show that for free N=2 VOSAs the construction contains as tensor factors both the free fermion vertex operator superalgebra and a parity-twisted module for the free fermion vertex operator, as studied, for instance in \cite{FFR}.  The  mirror-twisted modules for the other mirror map is related to permutation twisted constructions as developed by the author along with Dong and Mason in \cite{BDM} for VOAs, but extended to signed permutation automorphisms of tensor products of VOSAs. This extension of \cite{BDM} is nontrivial and has many interesting features and will be studied in detail in future  work.  

The construction and classification of Ramond twisted sectors for N=1 VOSAs includes or is related to results previously presented in several works, such as \cite{FFR}, \cite{Li-twisted}, \cite{S}, \cite{DZ}, \cite{Milas}.   Our construction of the $\sigma_\xi$-twisted N=2 VOSA modules uses results previously presented in \cite{Li-twisted}.  

We believe the current work should have interesting applications toward understanding and extending the work in, for instance, \cite{FFR}, \cite{Hohn-thesis}, \cite{Hohn-2008}, \cite{Duncan07}--\cite{Duncan-Ru}, to interesting N=2 supersymmetric settings. 

Much of this work is written in an introductory style, as several sections were first written as part of the lecture notes for the course ``Geometric and Algebraic Aspects of Superconformal Field Theory" taught by the author at the University of Notre Dame in Spring 2010.

Acknowledgment: 
The author gratefully thanks the NSA for the grant MSPF-07G-169 that partially supported this work.  In addition, the author thanks the Max Planck Institute for Mathematics in Bonn, Germany for a research grant that also partially supported this work, and thanks them for their hospitality during the 2010-2011 academic year.

\subsection{The N=1 Neveu-Schwarz algebra and the N=1 Ramond algebra}\label{N=1-NS-R-intro-section}

The {\it N=1 Neveu-Schwarz algebra} is the Lie superalgebra with basis consisting of the central element $d$, even elements $L_n$ for $n \in \mathbb{Z}$, and odd elements $G_{r}$ for $r \in \mathbb{Z} + \frac{1}{2}$, and supercommutation relations  
\begin{eqnarray}
\left[L_m ,L_n \right] &=& (m - n)L_{m + n} + \frac{1}{12} (m^3 - m) \delta_{m + n 
, 0} \; d , \label{Virasoro-relation-N1} \\
\left[ L_m, G_{r} \right] &=& \left(\frac{m}{2} - r\right) G_{m+r} ,\\
\left[ G_{r} , G_{s} \right] &=&2L_{r + s}  + \frac{1}{3} \left(r^2 - \frac{1}{4} \right) \delta_{r+s , 0} \; d ,  \label{N1-Neveu-Schwarz-relation-last}
\end{eqnarray}
for $m, n \in \mathbb{Z}$, and $r,s \in \mathbb{Z} + \frac{1}{2}$.  

The N=1 Ramond algebra is the Lie superalgebra with basis consisting of the central element $d$, even elements $L_n$ for $n \in \mathbb{Z}$,  and odd elements $G_r$ for $r \in \mathbb{Z}$, and supercommutation relations given by (\ref{Virasoro-relation-N1})--(\ref{N1-Neveu-Schwarz-relation-last}), where now $r, s \in \mathbb{Z}$.

Note that the only nontrivial Lie superalgebra automorphism of the N=1 Neveu-Schwarz algebra (resp. N=1 Ramond algebra), is the parity automorphism which is the identity on the even subspace (the Virasoro Lie algebra) and acts as $-1$ on the odd subspace (the subspace spanned by $G_r$ for either $r \in \mathbb{Z} + \frac{1}{2}$ in the Neveu-Schwarz case or $r \in \mathbb{Z}$ in the Ramond case).

Here we give some connections to the geometry underlying two-dimensional N=1 superconformal field theory and some motivation for our interest in the N=1 Neveu-Schwarz and Ramond algebras.  Let $k = 1$ or $\frac{1}{2}$, and let $x^k$ a commuting formal variable and $\varphi$ an anti-commuting formal variable.  The N=1 Neveu-Schwarz and N=1 Ramond algebras have the following representation with central element zero, in terms of superderivations on $\mathbb{C}[[x^k, x^{-k}]][\varphi]$, for $k =1$ and $\frac{1}{2}$, respectively:
\begin{eqnarray}
L_n(x,\varphi) &=&  - \biggl( x^{n+ 1} \frac{\partial}{\partial x} + \Bigl(\frac{n + 1}{2} \Bigr)x^n \varphi  \frac{\partial }{\partial \varphi }   \biggr) \label{L-notation-N1} \\
G_r (x,\varphi) &=&  - x^{r + 1/2} \Bigl(   \frac{\partial }{\partial \varphi}  -  \varphi  \frac{\partial}{\partial x} \Bigr)    \label{G-notation-N1}
\end{eqnarray} 
where $n \in \mathbb{Z}$, and $r \in \mathbb{Z} + \frac{1}{2}$ for the N=1 Neveu-Schwarz algebra, and $r \in \mathbb{Z}$ for the N=1 Ramond algebra.  The representation of the N=1 Ramond algebra in terms of superderivations is obtained from the representation of the N=1 Neveu-Schwarz algebra via the nonsuperconformal change of variables $(x, \varphi) \mapsto (x, \varphi x^{1/2})$ for the odd components. 

The superderivations $L_n(x,\varphi)$ and $G_r(x, \varphi)$ for $n \in \mathbb{Z}$ and $r \in \mathbb{Z} + \frac{1}{2}$ which give a representation of the N=1 Neveu-Schwarz algebra with central charge zero are the infinitesimal superderivations which give the data for genus-zero worldsheets on the supersphere corresponding to worldsheets with tubes swept out by a superstring propagating through space-time for two-dimensional holomorphic N=1 superconformal field theory where those tubes are anti-periodic in the fermionic components \cite{B-thesis}, \cite{B-memoirs}, \cite{B-iso-thm}.  To describe vertex operators on the supercylinder that is periodic in the fermionic components or to describe genus one and higher genus superstring interactions using the genus zero interactions, one needs to move to a twisted module of the underlying N=1 VOSA which is a representation of the N=1 Ramond algebra.   That is, extending the work of Zhu \cite{Zhu} to the supersymmetric setting, the N=1 superconformal change of variables needed to move from vertex operators on the superdisc to vertex operators on the anti-periodic supercylinder is $(z, \theta) \mapsto (e^z, e^{z/2} \theta)$ which is anti-periodic in the $\theta$ component, i.e. with periodicity given by $(z, \theta) \sim (z + 2 \pi i n, (-1)^n \theta)$.  The change of variables needed to move from vertex operators on the superdisc to vertex operators on the periodic supercylinder is the composition of the N=1 superconformal map from the superdisc to the supercylinder and the nonsuperconformal change of variables $(z, \theta) \mapsto (z, z^{1/2} \theta)$.  This composition gives the change of coordinates $(z, \theta) \mapsto (e^z, e^z \theta)$ from the a double cover of the disc to the supercylinder which is periodic in the $\theta $ component, i.e. with periodicity given by $(z, \theta) \sim (z + 2 \pi i n, \theta)$.  

This is one of the main motivations for studying the N=1 Neveu-Schwarz algebra and constructing certain twisted modules for N=1 supersymmetric VOSAs, since the appropriate twist gives rise to twisted modules which are representations of the N=1 Ramond algebra.  

\subsection{The N=2 superconformal algebras}\label{N=2-NS-R-intro-section}

The {\it N=2 Neveu-Schwarz Lie superalgebra} (also called the {\it N=2 superconformal algebra}) is the Lie superalgebra with basis consisting of the central element $d$, even elements $L_n$ and $J_n$ for $n \in \mathbb{Z}$, and odd elements $G^{(j)}_r$ for $j = 1,2$ and $r \in \mathbb{Z} + \frac{1}{2}$, and such that the supercommutation relations are given as follows:  $L_n$, $d$ and $G^{(j)}_r$ satisfy the supercommutation relations for the N=1 Neveu-Schwarz Lie superalgebra given by (\ref{Virasoro-relation-N1})--(\ref{N1-Neveu-Schwarz-relation-last}) for both $G_r = G^{(1)}_r$ and for $G_r = G^{(2)}_r$;  the remaining supercommutation relations are given by
\begin{eqnarray}
\left[L_m, J_n \right] &=& -n J_{m+n}, \qquad \quad \qquad \left[ J_m, J_n \right] \ =  \ \frac{1}{3} m \delta_{m+n,0} d \label{N2-evens}\\
\left[ J_m, G^{(1)}_{r} \right] \!  &=& \!  - i G^{(2)}_{m+r} , \qquad  \quad \quad  \, \left[ J_m, G^{(2)}_{r} \right]  \  = \   i G^{(1)}_{m+r} ,\label{nonhomo-J-relation} \\
\qquad \ \left[ G^{(1)}_r, G^{(2)}_s \right]  \! 
&=& \!  - i  (r-s) J_{r+s} . \label{N2-last} 
\end{eqnarray}

The {\it N=2 Ramond algebra} is the Lie superalgebra with basis consisting of the central element $d$, even elements $L_n$ and $J_n$ for $n \in \mathbb{Z}$, and odd elements $G_r^{(j)}$ for $r \in \mathbb{Z}$ and $j=1,2$, and supercommutation relations given by those of the N=2 Neveu-Schwarz algebra but with $r,s \in \mathbb{Z}$, instead of $r,s \in \mathbb{Z} + \frac{1}{2}$. 

More generally, there is an infinite family of algebras which includes the N=2 Neveu-Schwarz and the N=2 Ramond algebra and which are all isomorphic to the N=2 Neveu-Schwarz algebra under the so called ``spectral flow".  However it is easiest to express this phenomenon if we make a change of basis which is ubiquitous in superconformal field theory.
So consider the substitutions
\begin{equation}
G^{(1)}_r = \frac{1}{\sqrt{2}} \left( G^+_r + G^-_r \right), \qquad
G^{(2)}_r= \frac{i}{\sqrt{2}} \left( G^+_r - G^-_r \right) , 
\end{equation}
or equivalently  $G^\pm_r = \frac{1}{\sqrt{2}}(G^{(1)}_r \mp i G^{(2)}_r)$.  This substitution is equivalent to the change of variables $\varphi^\pm = \frac{1}{\sqrt{2}} ( \varphi^{(1)} \pm i \varphi^{(2)})$ in the variables $(x, \varphi^{(1)}, \varphi^{(2)})$ representing the one even and two odd local coordinates on an N=2 superconformal worldsheet representing superstrings propagating in space-time in N=2 superconformal field theory, see for instance \cite{B-n2moduli}.   The N=2 Neveu-Schwarz algebra or N=2 Ramond algebra is often written using these substitutions, so that the basis consists of the even central element $d$, the even elements $L_n$, $J_n$,  for $n \in \mathbb{Z}$, the odd elements $G^{\pm}_r$, for $r \in \mathbb{Z} + \frac{1}{2}$ for the N=2 Neveu-Schwarz algebra, or for $r\in\mathbb{Z}$ for the N=2 Ramond algebra, and supercommutation relations given by (\ref{Virasoro-relation-N1}), (\ref{N2-evens}) and 
\begin{eqnarray}
\left[ L_m, G^\pm_r \right] \! \! &=& \! \! \left(\frac{m}{2} -r \right) G^\pm_{m+r} , \label{Neveu-Schwarz-relation-first} \\
\left[ J_m, G^\pm_r \right] \! \! &=& \! \!  \pm G^\pm_{m+r} , \qquad \qquad \qquad \ \ \ \left[ G^\pm_r, G^\pm_s \right] \ =\  0 , \label{J-relation}\\
\qquad \left[ G^+_r , G^-_s\right]  \! \! &=& \! \! 2L_{r+ s} + (r-s) J_{r+s}  + \frac{1}{3} (r^2 -\frac{1}{4}) \delta_{r+s , 0} \; d ,\label{Neveu-Schwarz-relation-last} 
\end{eqnarray}
for $m, n \in \mathbb{Z}$,  and $r,s \in \mathbb{Z} + \frac{1}{2}$ for the N=2 Neveu-Schwarz algebra, or $r,s \in \mathbb{Z}$ for the N=2 Ramond algebra.  We call this the {\it homogeneous basis} for the N=2 Neveu-Schwarz or N=2 Ramond algebras. 

Observe that there is also the notion of a Lie superalgebra generated by even elements $L_n$ and $J_n$ for $n \in \mathbb{Z}$ and by odd elements $G^\pm_{r \pm t}$, for $r \in \mathbb{Z} + \frac{1}{2}$ and for $t \in \mathbb{C}$ (or in general any underlying field of characteristic zero which contains the rationals).  We shall call this algebra the {\it $t$-shifted N=2 superconformal algebra} or {\it $t$-shifted N=2 Neveu-Schwarz algebra}.  Thus the $t$-shifted N=2 Neveu-Schwarz algebra is the N=2 Neveu-Schwarz algebra if $t \in \mathbb{Z}$, and is the N=2 Ramond algebra if  $t \in \mathbb{Z} + \frac{1}{2}$.  As was first shown in \cite{SS}, the $t$-shifted N=2 Neveu-Schwarz algebras are all isomorphic under the continuous family of {\it spectral flow} maps, denoted $\mathcal{D}(t)$, for $t \in \mathbb{C}$, but which fix the Virasoro algebra only for $t = 0$.  These are given by 
\begin{equation}
\mathcal{D}(t) :  \quad \begin{array}{ll}
\displaystyle{L_n \ \mapsto \ L_n + tJ_n + \frac{t^2}{6}\delta_{n,o} d,  } &  \qquad d  \ \mapsto \ d,\\
\displaystyle{J_n \ \, \mapsto\  J_n + \frac{t}{3} \delta_{n,0} d, } & \quad G^{\pm}_r \ \mapsto \ G^{\pm}_{r \pm t}.   
\end{array}
\end{equation}

We shall show in Section \ref{other-iso-section} that for $t$ a positive rational number less than one,  representations of the $t$-shifted N=2 Neveu-Schwarz algebras  naturally occur as twisted modules of N=2 vertex operator superalgebras under twists arising from the group of automorphisms of the N=2 Neveu-Schwarz algebra which preserve the Virasoro algebra.

The group of automorphisms of the N=2 Neveu-Schwarz algebra (or more generally the $t$-shifted N=2 superconformal algebras)  which preserve the Lie subalgebra generated by $L_n$ and $J_n$ for $n \in \mathbb{Z}$ are given by:
\begin{equation}\label{sigmas}
\sigma_\xi:  \quad G^\pm_r \mapsto \xi^{\pm 1} G^\pm_r, \qquad J_n \mapsto J_n, \qquad L_n \mapsto L_n, \qquad d \mapsto d,
\end{equation}
for $\xi \in \mathbb{C}^\times$, if we are taking the algebra over $\mathbb{C}$, or more generally, for $\xi$ an invertible even element of the underlying base space.   In addition, we have the Virasoro-preserving automorphism which is commonly referred to as the {\it mirror map} given by:
\begin{equation}\label{mirror-map}
\kappa:   \quad G^\pm_r  \mapsto G^\mp_r, \qquad J_n \mapsto -J_n, \qquad L_n \mapsto L_n, \qquad d \mapsto d.
\end{equation}
The family $\sigma_\xi$ along with $\kappa$ generate all the Virasoro-preserving automorphisms of the N=2 Neveu-Schwarz algebra, and thus this group is isomorphic to $\mathbb{Z}_2 \times \mathbb{C}^\times$, cf. \cite{B-n2axiomatic}.  

In terms of the nonhomogeneous basis, these automorphisms are given by 
\begin{eqnarray}\label{nonhomo-auto3-hyperbolic}
G^{(1)}_r  \! \! &\mapsto& \! \!  (\cosh \beta) G^{(1)}_r + i (\sinh \beta) G^{(2)}_r, \nonumber \\
\sigma_\xi: \quad  G^{(2)}_r \! \! &\mapsto& \! \! i(\sinh \beta) G^{(1)}_r - (\cosh \beta) G^{(2)}_r , \\
J_n \! \! &\mapsto& \! \!  J_n, \qquad L_n  \mapsto  L_n, \qquad d  \mapsto d, \nonumber
\end{eqnarray}
for $\xi \in \mathbb{C}^\times$ and $e^\beta = \xi$, and by 
\begin{equation}\label{mirror-map-nonhomo}
\kappa: \quad  G^{(1)}_r \mapsto G^{(1)}_r, \qquad G^{(2)}_r \mapsto - G^{(2)}_r,  \qquad J_n \mapsto - J_n, \qquad  L_n \mapsto L_n,  \qquad d \mapsto d.
\end{equation}

In Sections \ref{ramond-twisted-section} and \ref{other-iso-section}, we show that given an N=2 vertex operator superalgebra $V$, then for $\eta = e^{2 \pi i /k}$ and $\xi = \eta^j$, if the automorphism $\sigma_\xi$ on the representation of the N=2 Neveu-Schwarz algebra extends to a vertex operator superalgebra automorphism of $V$ (which it always does when $\xi = -1$ or when $V$ admits a $J(0)$-grading by charge) then a $\sigma_\xi$-twisted $V$-module is a representation of the $\frac{j}{k}$-shifted N=2 Neveu-Schwarz algebra.   

We also show in Section \ref{mirror-twisted-section} that if the mirror automorphism $\kappa$ on the representation of the N=2 Neveu-Schwarz algebra on $V$ extends to a vertex operator superalgebra automorphism on $V$, then the $\kappa$-twisted $V$-module is a representation of the {\it mirror-twisted N=2 Neveu-Schwarz algebra}.  The {\it mirror-twisted N=2 Neveu-Schwarz algebra} is defined to be the Lie superalgebra with basis consisting of even elements $L_n$, and $J_r$ and central element $d$, odd elements $G^{(1)}_r$ and $G^{(2)}_n$, for $n \in \mathbb{Z}$ and $r \in \mathbb{Z} + \frac{1}{2}$, and supercommutation relations given as follows:  The $L_n$ and $G^{(1)}_r$ satisfy the supercommutation relations for the N=1 Neveu-Schwarz algebra with central charge $d$;  the $L_n$ and $G^{(2)}_n$ satisfy the supercommutation relations for the N=1 Ramond algebra with central charge $d$;  and the remaining supercommutation relations are 
\begin{eqnarray}\label{mirror-twisted-NS}
\left [L_n, J_r\right]  &=&  -r J_{n+r}, \qquad  \qquad \quad \ \left[ J_r, J_s \right]  \ \ =  \ \  \frac{1}{3} r \delta_{r+s,0} d\\
\left[ J_r, G^{(1)}_s \right] &=&  - i G^{(2)}_{r+s} ,  \qquad \quad \quad \, \left[ J_r, G^{(2)}_n \right]  \ \   =   \ \  i G^{(1)}_{r+n} ,\\
& & \hspace{-.1in} \left[ G^{(1)}_r, G^{(2)}_n \right]   \ \ = \ \  - i  (r-n) J_{r+n} . \label{mirror-twisted-NS-last}
\end{eqnarray}
Note that this mirror-twisted N=2 Neveu-Schwarz algebra is not isomorphic to the ordinary N=2 Neveu-Schwarz algebra \cite{SS}.

We tie this into the geometry underlying N=2 superconformal field theory, by recalling from, for instance \cite{B-n2moduli}, that formally, the infinitesimal N=2 superconformal transformations are given by the even superderivations in $\mathrm{Der}( \mathbb{C} [[x, x^{-1}]][\varphi^+, \varphi^-])$
\begin{eqnarray}
L_n(x,\varphi^+,\varphi^-) &=& - \biggl( x^{n + 1} \frac{\partial}{\partial x} + \Bigl(\frac{n + 1}{2}\Bigr) x^n \Bigl( \varphi^+ \frac{\partial}{\partial \varphi^+} + \varphi^- \frac{\partial}{\partial \varphi^-}\Bigr)\biggr) \label{L-notation}\\
J_n(x,\varphi^+,\varphi^-) &=& - x^n\Bigl(\varphi^+\frac{\partial}{\partial \varphi^+} - \varphi^- \frac{\partial}{\partial \varphi^-}\Bigr)  \label{J-notation}
\end{eqnarray}
and the odd superderivations
\begin{eqnarray}
G^\pm_{n -\frac{1}{2}} (x,\varphi^+,\varphi^-) = - \biggl( x^n \Bigl( \frac{\partial}{\partial \varphi^\pm} - \varphi^\mp \frac{\partial}{\partial x}\Bigr) \pm n x^{n-1} \varphi^+ \varphi^- \frac{\partial}{\partial \varphi^\pm} \biggr) \label{G-notation}
\end{eqnarray}
for $n \in \mathbb{Z}$.   These superderivations give a representation of the N=2 Neveu-Schwarz algebra with central charge zero.  Performing the coordinate transformation $(x, \varphi^+, \varphi^-) \mapsto (x, x^t \varphi^+, x^{-t} \varphi^-)$ on the odd components transforms the superderivations above to a representation of the $t$-shifted N=2 Neveu-Schwarz algebra, and $(x, \varphi^+, \varphi^-) \mapsto (x, \varphi^-, \varphi^+)$ corresponds to the mirror automorphism on the N=2 Neveu-Schwarz algebra representation.

The superderivations $L_n(x,\varphi^+, \varphi^-)$, $J_n(x,\varphi^+, \varphi^-)$ and $G^\pm_r(x, \varphi^+, \varphi^-)$ for $n \in \mathbb{Z}$ and $r \in \mathbb{Z} + \frac{1}{2}$ which give a representation of the N=2 Neveu-Schwarz algebra with central charge zero are the infinitesimal superderivations which give the data for genus-zero worldsheets on the N=2 supersphere corresponding to worldsheets with tubes swept out by a superstring propagating through space-time for two-dimensional holomorphic N=2 superconformal field theory where those tubes have certain fixed periodicity condition in the two fermionic components.  To describe vertex operators on the supercylinder that have different periodicity conditions in the fermionic components or to describe genus one and higher genus superstring interactions using the genus zero interactions, one needs to move to various twisted modules of the underlying N=2 VOSA.    This is in analogy with the N=1 case, but with a much higher degree of complexity due to the richer geometric structure an N=2 superconformal surfaces possess, cf.  \cite{B-uniformization}, \cite{B-autogroups}.  

\subsection{Summary of results}
We summarize some of the main results of this paper as follows.  In Sections \ref{ramond-twisted-section}--\ref{other-iso-section}, we show:

\begin{thm}\label{first-thm}
If the Virasoro-preserving automorphisms of the N=2 Neveu-Schwarz algebra $\kappa$, and $\sigma_\xi$, for $\xi$ a root of unity, extend to VOSA automorphisms for an N=2 VOSA, $V$, then:\\
(i) A weak $\kappa$-twisted $V$-module is a representation of the mirror-twisted N=2 superconformal algebra.\\
(ii) A weak $\sigma_\xi$-twisted $V$-module, for $\xi = e^{2j \pi i /k}$, is a representation of the $\frac{j}{k}$-shifted N=2 superconformal algebra.  If $\xi = -1$, then such a VOSA automorphism always exists (the parity map), and in this case, a weak $\sigma_{-1}$-twisted $V$-module is a representation of the N=2 Ramond algebra.  
\end{thm}

In Sections \ref{ramond-construction-section}--\ref{sigma-twisted-section}, we show in particular:

\begin{thm}\label{second-thm}
 If $V$ is a free or lattice N=2 VOSA, then each Virasoro-preserving automorphism of the N=2 Neveu-Schwarz algebra extends to a VOSA automorphism of $V$, but not uniquely in the case of the mirror map.  

For free and lattice N=2 VOSAs, there are two distinct mirror maps representing two distinct conjugacy classes in the group of automorphisms of the VOSA, and thus giving rise to non-isomorphic mirror-twisted $V$-module structures.
\end{thm}

\subsection{The notions of vertex operator superalgebra, and N=1 or N=2 supersymmetric vertex operator superalgebra}\label{vosa-definitions-section}

In this section, we recall the notion of vertex operator superalgebra, as well as N=1 or N=2 Neveu-Schwarz vertex operator superalgebra, following
the notation and terminology of \cite{B-vosas}, \cite{B-iso-thm} and \cite{B-n2axiomatic}. Let $x, x_0, x_1, x_2,$ etc., denote commuting independent formal variables.
Let $\delta (x) = \sum_{n \in \Z} x^n$.  We will use the binomial
expansion convention, namely, that any expression such as $(x_1 -
x_2)^n$ for $n \in \C$ is to be expanded as a formal power series in
nonnegative integral powers of the second variable, in this case
$x_2$.

A {\it vertex operator superalgebra} is a $\frac{1}{2}\Z$-graded vector space
\begin{equation}
V=\coprod_{n\in \frac{1}{2} \Z}V_n
\end{equation}
satisfying ${\rm dim} \, V < \infty$ and $V_n = 0$ for $n$ sufficiently negative, 
that is also $\Z_2$-graded by {\it sign} 
\[V = V^{(0)} \oplus V^{(1)},\]
and equipped with a linear map
\begin{eqnarray}
V &\longrightarrow& (\mbox{End}\,V)[[x,x^{-1}]]\\
v &\mapsto& Y(v,x)= \sum_{n\in\Z}v_nx^{-n-1} \nonumber
\end{eqnarray}
and with two distinguished vectors ${\bf 1}\in V_0$, (the {\it vacuum vector})
and $\omega\in V_2$
(the {\it conformal element}) satisfying the following conditions for $u, v \in
V$:
\begin{eqnarray}
& &u_nv=0\ \ \ \ \ \mbox{for $n$ sufficiently large};\\
& &Y({\bf 1},x)=1;\\
& &Y(v,x){\bf 1}\in V[[x]]\ \ \ \mbox{and}\ \ \ \lim_{x\to
0}Y(v,x){\bf 1}=v;
\end{eqnarray}
\begin{multline}
x^{-1}_0\delta\left(\frac{x_1-x_2}{x_0}\right) Y(u,x_1)Y(v,x_2) - (-1)^{|u||v|}
x^{-1}_0\delta\left(\frac{x_2-x_1}{-x_0}\right) Y(v,x_2)Y(u,x_1)\\
= x_2^{-1}\delta \left(\frac{x_1-x_0}{x_2}\right) Y(Y(u,x_0)v,x_2)
\end{multline}
(the {\it Jacobi identity}), where $|v| = j$ if $v \in V^{(j)}$ for $j \in \Z_2$;
\begin{equation}
[L(m),L(n)]=(m-n)L(m+n)+\frac{1}{12}(m^3-m)\delta_{m+n,0}c
\end{equation}
for $m, n\in \Z,$ where
\begin{equation}
L(n)=\omega_{n+1}\ \ \ \mbox{for $n\in \Z$, \ \ \ \ i.e.},\
Y(\omega,x)=\sum_{n\in\Z}L(n)x^{-n-2}
\end{equation}
and $c \in \mathbb{C}$ (the {\it central charge} of $V$);
\begin{eqnarray}
& &L(0)v=nv=(\mbox{wt}\,v)v \ \ \ \mbox{for $n \in  \frac{1}{2} \Z$ and $v\in V_n$}; \\
& &\frac{d}{dx}Y(v,x)=Y(L(-1)v,x). 
\end{eqnarray}
This completes the definition. We denote the vertex operator superalgebra (VOSA) just
defined by $(V,Y,{\bf 1},\omega)$, or briefly, by $V$.

As a consequence of the definition, we have that 
\begin{equation}\label{consequences0}
\omega = L(-2) \mathbf{1} \quad \mbox{and} \quad L(n) \mathbf{1} = 0 \quad \mbox{for $n \geq -1$},
\end{equation}  
as well as 
\begin{equation}\label{L(-1)-action}
L(-1)v = v_{-2} \mathbf{1}.
\end{equation}


If a vertex operator superalgebra, $(V, Y, {\bf 1}, \omega)$, contains an element $\tau \in V_{3/2}$ satisfying 
\[Y(\tau, z) = \sum_{n \in \Z} \tau_n x^{-n-1} = \sum_{n \in \Z} G(n + 1/2) x^{-n-2} ,\]
where the $G(n + 1/2) = \tau_{n+1} \in (\mathrm{End} (V))^{(1)}$ generate a representation of the 
N=1 Neveu-Schwarz Lie superalgebra (that is $\frac{1}{2} G(-1/2)\tau = \omega$ with $L(n) = \omega_{n+1} \in (\mathrm{End} (V))^{(0)}$ which, along with the $G(n + 1/2)$ satisfy the N=1 Neveu-Schwarz Lie superalgebra relations (\ref{Virasoro-relation-N1})--(\ref{N1-Neveu-Schwarz-relation-last})), then we call $(V, Y, {\bf 1}, \tau)$ an {\it N=1 Neveu-Schwarz vertex operator superalgebra}, or a {\it N=1 supersymmetric vertex operator superalgebra}, or just N=1 VOSA for short. 

For an N=1 VOSA, it follows from the definition that 
\begin{equation}\label{consequences1}
\tau = G(-3/2) \mathbf{1} \quad \mbox{and} \quad G(n + 1/2) \mathbf{1} = 0 \quad \mbox{for $n \geq -1$.}
\end{equation}

If a VOSA $(V, Y, \mathbf{1}, \omega)$ has two vectors $\tau^{(1)}$ and $\tau^{(2)}$ such that $(V, Y, \mathbf{1}, \tau^{(j)})$ is an N=1 VOSA for both $j=1$ and $j = 2$, and the $\tau^{(j)}_{n+1} = G^{(j)}(n + 1/2)$ generate a representation of the N=2 Neveu-Schwarz Lie superalgebra, then we call such a VOSA an {\it N=2 Neveu-Schwarz vertex operator superalgebra} or {\it N=2 supersymmetric vertex operator superalgebra}, or for short, an N=2 VOSA.   In particular, we have that if $V$ is an N=2 VOSA, then there exists a vector $\mu = \frac{i}{2} G^{(1)}(1/2) \tau^{(2)} = -\frac{i}{2} G^{(2)} (1/2) \tau^{(1)} \in V_{(1)}$ such that writing 
\begin{equation}
Y(\mu, x) = \sum_{n \in \Z} \mu_n x^{-n-1} = \sum_{n \in \Z} J(n) x^{-n-1}
\end{equation}
we have that the $J(n) \in (\mathrm{End} \, V)^0$ along with the $G^{(j)}(n+ 1/2)$ and $L(n) = \omega_{n+1}$ for $\omega = \frac{1}{2} G^{(j)}(-1/2) \tau^{(j)}$ satisfy the supercommutation relations for the N=2 Neveu-Schwarz Lie superalgebra.

For an N=2 NS-VOSA, it follows from the definition that 
\begin{equation}\label{consequences2}
\mu = J(-1) \mathbf{1} \quad \mbox{and} \quad J(n) \mathbf{1} = 0 \quad \mbox{for $n \geq 0$.}
\end{equation}

If $V$ is an N=2 VOSA such that $V$ is not only $\frac{1}{2} \mathbb{Z}$ graded by $L(0)$ but also $\mathbb{Z}$-graded by $J(0)$ such that $J(0)v = nv$ with $n \equiv j \,  \mathrm{mod} \, 2$, for $v \in V^{(j)}$ for $j = 0,1$, then we say that $V$ is {\it $J(0)$-graded} or {\it graded by charge}.   

Given two vertex operator superalgebras $(V_1, Y_1, \mathbf{1}^{(1)}, \omega^{(1)})$ and $(V_2, Y_2, \mathbf{1}^{(2)}, \omega^{(2)})$, we have that $(V_1 \otimes V_2, \, Y,  \, \mathbf{1}^{(1)} \otimes \mathbf{1}^{(2)}, \, \omega^{(1)} \otimes \mathbf{1}^{(2)} + \mathbf{1}^{(1)} \otimes \omega^{(2)})$ is a vertex operator superalgebra, where $Y$ is given by 
\begin{equation}
Y(u_1 \otimes u_2, x) (v_1 \otimes v_2) = (-1)^{|u_2||v_1|} Y_1(u_1, x) v_1\otimes Y_2(u_2, x)v_2,
\end{equation}
for $u_1 \otimes u_2, \, v_1 \otimes v_2 \in V_1 \otimes V_2$.  

\section{Twisted modules for N=1 and N=2 VOSAs}

In this section, we present the notion of twisted modules for a VOSA and determine the structure of the $g$-twisted modules for an N=2 VOSA when $g$ arises from a finite, Virasoro-preserving  automorphism of the N=2 Neveu-Schwarz algebra.

\subsection{Automorphisms of VOSAs and the notion of twisted VOSA-module}\label{twisted-definitions-section}

In this section, we recall the notion of a $g$-twisted $V$-module for a vertex operator superalgebra $V$ and an automorphism $g$ of $V$ of finite order following the notation of, for instance,  \cite{Li-twisted}, \cite{DLM1}, \cite{DLM2}, \cite{BDM}, \cite{BHL}.

An {\it automorphism} of a vertex operator superalgebra $V$ is a linear
automorphism $g$ of $V$ preserving ${\bf 1}$ and $\omega$ such that
the actions of $g$ and $Y(v,x)$ on $V$ are compatible in the sense
that
\begin{equation}\label{automorphism}
g Y(v,x) g^{-1}=Y(gv,x)
\end{equation}
for $v\in V.$ Then $g V_n\subset V_n$ for $n\in  \frac{1}{2} \mathbb{Z}$.

If $g$
has finite order, $V$ is a direct sum of the eigenspaces $V^j$ of $g$,
\begin{equation}
V=\coprod_{j\in \Z /k \Z }V^j,
\end{equation}
where $k \in \mathbb{Z}_+$ is a period of $g$ (i.e., $g^k = 1$ but $k$
is not necessarily the order of $g$) and
\begin{equation}
V^j=\{v\in V \; | \; g v= \eta^j v\},
\end{equation}
for $\eta$ a fixed primitive $k$-th root of unity.

Note that we have the following $\delta$-function identity
\begin{equation}\label{delta-identity}
x_2^{-1} \delta \left (\frac{x_1 - x_0}{x_2} \right) \left( \frac{x_1 - x_0}{x_2} \right)^{k} = x_1^{-1} \delta \left (\frac{x_2 + x_0}{x_1} \right) \left( \frac{x_2 + x_0}{x_1} \right)^{-k}
\end{equation}
for any $k \in \mathbb{C}$.

We next review the notions of weak, weak admissible and
ordinary $g$-twisted module for a vertex operator superalgebra $V$
and an automorphism $g$ of $V$ of finite order $k$. 

Let $(V,Y,{\bf 1},\omega)$ be a vertex operator superalgebra
and let $g$ be an automorphism of $V$ of period $k \in
\mathbb{Z}_+$. A {\it weak $g$-twisted $V$-module} is a 
vector space $M$ equipped with a linear map
\begin{eqnarray}
V &\longrightarrow & (\mbox{End}\,M)[[x^{1/k},x^{-1/k}]] \\
v &\mapsto & Y^g (v,x)=\sum_{n\in \frac{1}{k}\Z }v_n^g x^{-n-1} \nonumber
\end{eqnarray}
satisfying the following conditions for $u,v\in V$ and $w\in M$:
\begin{eqnarray}
\label{cosetrelation}
& &Y^g(v,x)=\sum_{n\in  \Z + \frac{j}{k}}v_n^g x^{-n-1}\ \ \ \
\mbox{for $j\in \Z/k\Z$ and $v\in V^j$}; \\
& &v_n^g w=0\ \ \ \mbox{for $n$ sufficiently large};\\
& &Y^g ({\bf 1},x)=1;
\end{eqnarray}
\begin{multline}\label{twisted-jacobi}
x^{-1}_0\delta\left(\frac{x_1-x_2}{x_0}\right)
Y^g(u,x_1)Y^g(v,x_2)-  (-1)^{|u||v|} x^{-1}_0\delta\left(\frac{x_2-x_1}{-x_0}\right) Y^g
(v,x_2)Y^g (u,x_1)\\
= x_2^{-1}\frac{1}{k}\sum_{j\in \Z /k \Z}
\delta\left(\eta^j\frac{(x_1-x_0)^{1/k}}{x_2^{1/k}}\right)Y^g (Y(g^j
u,x_0)v,x_2)
\end{multline}
(the {\it twisted Jacobi identity}) where $\eta$ is a fixed primitive
$k$-th root of unity.  

We denote a weak $g$-twisted $V$-module by $(M, Y^g)$, or briefly, by $M$.  

If we take $g=1$, then we obtain the notion of weak $V$-module.  Note that the notion of weak $g$-twisted $V$-module for a vertex operator superalgebra is equivalent to the notion of $g$-twisted $V$-module for $V$ as a vertex superalgebra, cf. \cite{Li-twisted}.  In particular, the term ``weak" simply implies that we are making no assumptions about a grading on $M$.

Formula (\ref{cosetrelation}) can be expressed as follows: For $v \in
V$,
\begin{equation}
Y^g (gv,x) = \lim_{x^{1/k} \rightarrow \eta^{-1} x^{1/k}} Y^g(v, x),
\end{equation}
where the limit stands for formal substitution. 

As a consequence of the definition, we have the following supercommutator relation on $M$ for $u \in V^j$:
\begin{multline}\label{super-commutator}
[Y^g(u,x_1), Y^g(v, x_2) ] \\
= \mathrm{Res}_{x_0} x_2^{-1} \delta\left(\frac{x_1-x_0}{x_2}\right)\left(\frac{x_1 - x_0}{x_2}\right)^{-j/k} Y^g (Y( u,x_0)v,x_2),
\end{multline}
which follows from taking $\mathrm{Res}_{x_0}$ of both sides of the twisted Jacobi identity (\ref{twisted-jacobi}).
In addition, multiplying both sides of (\ref{twisted-jacobi}) by $\left(\frac{x_1 - x_0}{x_2} \right)^{j/k}$, taking $\mathrm{Res}_{x_1}$ of both sides, and using the $\delta$-function identity (\ref{delta-identity}), we have the following formula for iterates for the $g$-twisted vertex operators on $M$ for $u \in V^j$:
\begin{multline}\label{super-iterate}
Y^g (Y( u,x_0)v,x_2) = \mathrm{Res}_{x_1} \left(\frac{x_1 - x_0}{x_2} \right)^{j/k} \left( x^{-1}_0\delta\left(\frac{x_1-x_2}{x_0}\right)
Y^g(u,x_1)Y^g(v,x_2) \right. \\
\left. -  (-1)^{|u||v|} x^{-1}_0\delta\left(\frac{x_2-x_1}{-x_0}\right) Y^g
(v,x_2)Y^g (u,x_1) \right).
\end{multline}

Letting $u = v = \omega$ and taking $\mathrm{Res}_{x_1} x_1^{m+1}  \mathrm{Res}_{x_2} x_2^{n+1}$ of both sides of the supercommutator relation (\ref{super-commutator}), and then using (\ref{consequences0}), the Virasoro relations for $L(n) \in \mathrm{End} \, V$ for $n \in \mathbb{Z}$, and the delta function identity (\ref{delta-identity}), it follows that for a weak $g$-twisted $V$-module, $M$, we have
\begin{equation}\label{Viralgrelations}
[L^g(m),L^g(n)]=(m-n)L^g (m+n)+\frac{1}{12}(m^3-m)\delta_{m+n,0}c
\end{equation}
for $m, n\in \Z$, where $c$ is the central charge of $V$, and
\begin{equation}\label{define-L^g}
L^g (n)=\omega_{n+1}^g \ \ \ \mbox{for $n\in \Z$, \ \ \ \ \ i.e.},\ Y^g
(\omega,x)=\sum_{n\in \Z} L^g (n)x^{-n-2}.
\end{equation}
In addition, letting $v = \mathbf{1}$, taking $\mathrm{Res}_{x_0} x_0^{-2}$ of both sides of the iterate formula (\ref{super-iterate}), using (\ref{L(-1)-action}), and an argument analogous to that in the proof of Prop. 3.2.18 in \cite{LL}, we have
\begin{equation}\label{twisted-L-derivative}
\frac{d}{dx}Y^g (u,x)=Y^g (L(-1)u,x).
\end{equation}

Let $(M_1, Y^g_1)$ and $(M_2, Y^g_2)$ be two weak $g$-twisted $V$-modules.  A {\it $g$-twisted $V$-module homomorphism} from $M_1$ to $M_2$, is a linear map $f: M_1 \longrightarrow M_2$ such that 
\begin{equation}
f(Y_1^g(v,x)w) = Y^g_2(v,x) f(w)
\end{equation}
for $v \in V$ and $w \in M_1$.

A {\em weak admissible} $g$-twisted $V$-module is a weak $g$-twisted
$V$-module $M$ which carries a $\frac{1}{2k}{\Z}$-grading
\begin{equation}\label{m3.12}
M=\coprod_{n\in\frac{1}{2k}\Z}M(n)
\end{equation}
such that $v^g_mM(n)\subseteq M(n+\mathrm{wt} \; v-m-1)$ for homogeneous $v\in V$, and $M(n) = 0$ for $n$ sufficiently small. 
If $g=1,$ we have the notion of weak admissible $V$-module.

\begin{rema}
{\em Above we used the term ``weak admissible $g$-twisted module''
whereas in much of the literature (cf. \cite{DLM1}, \cite{BDM}) the term
``admissible $g$-twisted module'' is used for this notion.  We used
the qualifier ``weak'' to stress that these are indeed only weak
modules and in general are not ordinary modules.  However, for the
sake of brevity, we will now drop the qualifier ``weak''.}
\end{rema}

The vertex operator superalgebra $V$ is called $g$-{\em rational} if every
admissible $g$-twisted $V$-module is completely reducible, i.e., a
direct sum of irreducible admissible $g$-twisted modules.

An {\it ordinary $g$-twisted $V$-module}  is a weak $g$-twisted $V$-module $M$ which is $\C$-graded
\begin{equation}
M=\coprod_{\lambda \in \C}M_\lambda
\end{equation}
such that for each $\lambda$, $\dim M_{\lambda}< \infty $ and $M_{n/k
+\lambda}=0$ for all sufficiently negative integers $n$.  In addition,
\begin{equation}\label{L^g-grading}
L^g (0) w=\lambda w \qquad \mbox{for $w \in M_\lambda$}.
\end{equation}
We will usually refer to an ordinary $g$-twisted $V$-module, as just a $g$-twisted $V$-module.  
We call a $g$-twisted $V$-module $M$ {\it simple} or {\it irreducible}
if the only submodules are 0 and $M$.  

For a $g$-twisted $V$-module, $M$, we have the notion of {\it graded dimension} or {\it $q$-dimension}, denoted $\mathrm{dim}_q M$, and defined to be
\begin{equation}
\mathrm{dim}_q M = tr_M q^{L^g(0) - c/24} = q^{-c/24} \sum_{\lambda \in \mathbb{C}} (\mathrm{dim} \, M_\lambda) q^\lambda.
\end{equation}

If $V$ is an N=2 VOSA, and $M$ is a $g$-twisted $V$-module such that each $M_\lambda$ is also $J(0)$-graded, then we also have the notion of $J(0)$- and $L(0)$-graded dimension, or 
{\it $p,q$-dimension} given by
\begin{equation}
\mathrm{dim}_{p,q} M = tr_M p^{J^g(0)} q^{L^g(0) - c/24} .
\end{equation}

\subsection{Twisting by the parity involution $\sigma$---the Ramond sectors}\label{ramond-twisted-section}

Let $V$ be a VOSA and $g = \sigma$ the parity automorphism of $V$ given by
\begin{eqnarray}
\sigma : V & \longrightarrow & V \\
v &\mapsto&  (-1)^{|v|} v.\nonumber
\end{eqnarray}
Then the eigenspaces of $\sigma$ are $V^j = V^{(j)}$ for eigenvalue $\eta^j = (-1)^j$, for $j \in \Z_2$.

Let $M$ be a weak $\sigma$-twisted $V$-module.  Now suppose that in addition to being a vertex operator superalgebra, $V$ is an N=1 Neveu-Schwarz vertex operator superalgebra with N=1 superconformal element $\tau \in V^1$.  Then from (\ref{automorphism}), we have that $\sigma v_n \sigma^{-1} = (-1)^j v_n$ and thus under $\sigma$ we have $L(n) \mapsto L(n)$ and $G(r) \mapsto -G(r)$, for $n \in \mathbb{Z}$ and $r \in \mathbb{Z} + \frac{1}{2}$.   That is, the VOSA automorphism $\sigma$ acts via conjugation as the parity Virasoro-preserving automorphism on the representation of the N=1 Neveu-Schwarz algebra realized by $V$.  Note that by a slight abuse of notation, we will refer to both of these maps as ``$\sigma$".

Write 
\begin{equation}\label{Ramond-twisted-operators}
Y^\sigma(\omega, x) = \sum_{n \in \mathbb{Z}}  L^\sigma(n) x^{-n-2}, \qquad \qquad Y^\sigma (\tau, x) = \sum_{n \in \mathbb{Z}} G^\sigma(n) x^{-n - 3/2}, 
\end{equation}
i.e., define $G^\sigma(n) \in \mathrm{End} (M)$, for $n \in \mathbb{Z}$, by $\tau^\sigma_{n + 1/2} = G^\sigma(n)$.
Then using the supercommutator relations (\ref{super-commutator}) for the twisted vertex operators acting on $M$ for $u = \tau$ and $v = \tau$ or $\omega$, using the $L (-1)$-derivative property for the twisted vertex operators (\ref{twisted-L-derivative}), using the N=1 Neveu-Schwarz supercommutation relations for $L(n), G(n + 1/2) \in \mathrm{End} (V)$, using the fact that $L(n) \mathbf{1} = G(n + 1/2) \mathbf{1} = 0$ for $n \geq -1$, and using the $\delta$-function identity (\ref{delta-identity}),  
we have that the odd endomorphisms $G^\sigma(n)$ and the even endomorphisms $L^\sigma(n)$ on $M$, for $n \in \mathbb{Z}$, satisfy the supercommutation relations for the N=1 Ramond algebra (\ref{Virasoro-relation-N1})--(\ref{N1-Neveu-Schwarz-relation-last}) with central charge $c$.

Similarly, if $V$ is an N=2 VOSA with superconformal elements $\tau^{(j)}$ for $j = 1,2$, and $M$ is a $\sigma$-twisted module for $V$, then we have that writing
\begin{equation}\label{N2-Ramond-twisted-operators}
\begin{array}{cc}
Y^\sigma(\omega, x) = \displaystyle{\sum_{n \in \mathbb{Z}}  L^\sigma(n) x^{-n-2}}, \qquad Y^\sigma(\mu, x) = \displaystyle{\sum_{n \in \mathbb{Z}} J^\sigma(n) x^{-n-1}}, \\
Y^\sigma (\tau^{(j)}, x) = \displaystyle{\sum_{n \in \mathbb{Z}} G^{(j), \sigma} (n) x^{-n - 3/2}}, 
\end{array}
\end{equation}
for $j =1,2$, the $G^{(j), \sigma}(n) = (\tau^{(j)}_{n + 1/2})^\sigma$, along with $L^\sigma(n) = \omega^\sigma_{n+1}$ and $J^\sigma(n) = \mu^\sigma_n$, for $n \in \Z$, generate a representation of the N=2 Ramond superalgebra with central charge $c$.  

Note that as a Virasoro-preserving automorphism of the N=2 Neveu-Schwarz algebra, $\sigma = \sigma_{-1}$ in the notation of (\ref{sigmas}).

\subsection{Mirror maps and mirror-twisted modules for N=2 VOSAs}\label{mirror-twisted-section}  

Recall from Section \ref{N=2-NS-R-intro-section}  that the N=2 Neveu-Schwarz Lie superalgebra has an automorphism $\kappa$ given by (\ref{mirror-map}) in the homogeneous basis or by (\ref{mirror-map-nonhomo}) in the nonhomogeneous basis.  

If an N=2 VOSA, $V$ with central charge $c$, has a VOSA automorphism, which we also call $\kappa$, such that $\kappa(\mu) = -\mu$, and $\kappa(\tau^\pm) = \tau^\mp$ (or equivalently $\kappa(\tau^{(1)}) = \tau^{(1)}$  and $\kappa(\tau^{(2)}) = -\tau^{(2)}$), then such a VOSA automorphism of $V$ is called an {\it N=2 VOSA mirror map}.   If such a map exists for $V$, then $\kappa$ acting by conjugation on $\mathrm{End} \, V$ restricts to the mirror map Virasoro-preserving automorphism $\kappa$ on the elements $L(n)$, $J(n)$, and $G^\pm (r)$, for $n \in \mathbb{Z}$ and $r \in \mathbb{Z} + \frac{1}{2}$, which give the N=2 Neveu-Schwarz algebra representation on $V$ generated by $\tau^\pm$.

In this case, if $M$ is a  weak $\kappa$-twisted module for $V$, then write $Y^{\kappa}$ for the $\kappa$-twisted operators, and
\begin{equation}\label{mirror-twisted-operators}
\begin{array}{llllll}
Y^\kappa(\omega, x) \! \! \! & = \displaystyle{\sum_{n \in \Z} L^\kappa (n) x^{-n-2}, } & Y^\kappa (\tau^{(1)}, x)  \! \! \!   & = \displaystyle{\sum_{r \in \Z + \frac{1}{2} } G^{(1),\kappa}(r) x^{-r-\frac{3}{2}} }\\
Y^\kappa(\mu, x)  \! \! \!   &=   \!  \displaystyle{\sum_{r \in \Z + \frac{1}{2} } J^\kappa (r) x^{-r-1} }, \ \ \ &  Y^\kappa (\tau^{(2)}, x)  \! \! \!  &=     \displaystyle{\sum_{n \in \Z} G^{(2),\kappa}(n) x^{-n-\frac{3}{2}} }.
\end{array}
\end{equation} 
That is, define $J^\kappa (n) = \mu_n^\kappa$ and $G^{(2), \kappa}(n - 1/2) = \tau^{(2),\kappa}_n$, for $n \in \mathbb{Z} + \frac{1}{2}$.  Then,
using the supercommutator relations (\ref{super-commutator}) for the $\kappa$-twisted vertex operators acting on $M$, using the $L(-1)$-derivative property and the N=2 Neveu-Schwarz supercommutation relations, using the fact that $L(n) \mathbf{1} = G^{(j)}(n + 1/2) \mathbf{1} = J(n+1) \mathbf{1} = 0$ for $n \geq -1$ and for $j = 1,2$, we have that the supercommutation relations for the $\kappa$-twisted modes of $\omega$, $\mu$, $\tau^{(1)}$ and $\tau^{(2)}$, given by $L^\kappa(n), G^{(2), \kappa}(n)$, for $n \in \mathbb{Z}$, and $J^\kappa(r), G^{(1), \kappa}(r)$, for $r \in \mathbb{Z} + \frac{1}{2}$,  satisfy the relations of the mirror-twisted N=2 Neveu-Schwarz algebra given by (\ref{mirror-twisted-NS})--(\ref{mirror-twisted-NS-last}) with central charge $c$.

In particular, a $\kappa$-twisted module, $M$, for an N=2 VOSA reduces the N=2 Neveu-Schwarz algebra representation to an N=1 Neveu-Schwarz algebra representation coupled with an N=1 Ramond algebra representation.

\subsection{Twisting by an automorphism corresponding to the N=2 Neveu-Schwarz algebra automorphism $\sigma_\xi$}\label{other-iso-section}

In addition to the parity map and the mirror map, the N=2 Neveu-Schwarz Lie superalgebra has the automorphisms $\sigma_\xi$ given by (\ref{sigmas}) for $\xi \in \mathbb{C}^\times$ and $\xi \neq -1$.   Of course if $\xi = -1$, this is just the parity map; but for $\xi \neq \pm1$, this gives additional automorphisms of the N=2 Neveu-Schwarz algebra which can possibly be extended to VOSA automorphisms of an N=2 VOSA, $V$.    If such an automorphism is of finite order, i.e. if $\xi$ is a root of unity, and $\sigma_\xi$ extends to a VOSA automorphism of $V$, then we have the notion of a $\sigma_\xi$-twisted $V$-module.  

For instance, if $V$ is an N=2 VOSA which is also $J(0)$-graded such that the $J(0)$ eigenvalues are integral with $J(0) \omega = J(0) \mu = 0$, and $J(\tau^{(\pm)}) = \pm \tau^{(\pm)}$, then setting $\sigma_\xi(v) = \xi^n v$ if $J(0)v = nv$ gives such a VOSA automorphism, and this automorphism, acting by conjugation on $\mathrm{End} \, V$ restricts to the Virasoro-preserving automorphism, also denoted $\sigma_\xi$, on the representation of the N=2 Neveu-Schwarz algebra.

Let $\eta = e^{2\pi i /k}$, for $k \in \mathbb{Z}_+$, and let $\xi = \eta^j$, for $j = 1, \dots, k-1$.  Let $\sigma_\xi$ be a VOSA automorphism of an N=2 VOSA, $V$, such that $\sigma_\xi (\mu) = \mu$ and $\sigma_\xi(\tau^{(\pm)}) = \xi^{\pm1} \tau^{(\pm)}$.  Then $\omega, \mu \in V^{0}$ and $\tau^{(\pm)} \in V^{\pm j}$.  
 If such a map exists for $V$, and $M$ is a  weak $\sigma_\xi$-twisted module for $V$, then write $Y^{\sigma_\xi}$ for the $\sigma_\xi$-twisted operators, and
\begin{equation}\label{sigma-twisted-operators}
\begin{array}{cc}
Y^{\sigma_\xi}(\omega, x)  = \displaystyle{\sum_{n \in \Z} L^{\sigma_\xi} (n) x^{-n-2}, } \qquad \qquad Y^{\sigma_\xi} (\mu, x)  = \displaystyle{\sum_{n \in \Z} J^{\sigma_\xi}(n) x^{-n-1} }\\
Y^{\sigma_\xi}(\tau^{(\pm)}, x)  =   \displaystyle{\! \! \! \sum_{r \in \Z - \frac{1}{2} \pm  \frac{j}{k} } \! \! \! \! G^{\pm, \sigma_\xi} (r) x^{-r-\frac{3}{2}} }. 
\end{array}
\end{equation} 
Then using the supercommutator relations  (\ref{super-commutator})  for the $\sigma_\xi$-twisted vertex operators acting on $M$, using the $L(-1)$-derivative property and the N=2 Neveu-Schwarz supercommutation relations on $V$, and using the fact that $L(n) \mathbf{1} = G^{\pm}(n + 1/2) \mathbf{1} = J(n+1) \mathbf{1} = 0$ for $n \geq -1$, we have that the supercommutation relations for the $\sigma_\xi$-twisted modes of $\omega$, $\mu$, and  $\tau^{(\pm)}$, that is the $L^{\sigma_\xi}(n)$ and $J^{\sigma_\xi}(n)$ for $n \in \mathbb{Z}$, and $G^\pm(r)$ for $r \in \mathbb{Z} + \frac{1}{2} \pm \frac{j}{k}$, respectively, satisfy the relations for the $\frac{j}{k}$-shifted N=2 Neveu-Schwarz algebra (\ref{Virasoro-relation-N1}), (\ref{N2-evens}), (\ref{Neveu-Schwarz-relation-first})--(\ref{Neveu-Schwarz-relation-last}) with central charge $c$. 

That is the sectors for N=2 supersymmetric VOSAs that arise under fractional spectral flow $\mathcal{D}(t)$, for $t = j/k$, $k \in \mathbb{Z}_+$, $j = 1, \dots, k-1$,  are twisted sectors under the Virasoro-preserving automorphisms $\sigma_\xi$ of the N=2 Neveu-Schwarz algebra.

\section{Free and lattice constructions of N=1 and N=2 VOSAs}\label{free-section}

We follow the notation of \cite{LL} for free and lattice VOAs and give explicitly the elementary constructions of free bosonic, free fermionic and free N=1 and N=2 VOSAs.  

Let $\mathfrak{h} = \mathfrak{h}^0 \oplus \mathfrak{h}^1$ be a Lie superalgebra.   If  
\[ \mathrm{Cent} \, \mathfrak{h} = [\mathfrak{h}, \mathfrak{h}], \quad \mathrm{and} \quad \mathrm{dim \, Cent} \, \mathfrak{h} = 1,\]
then $\mathfrak{h}$ is said to be a {\it Heisenberg Lie superalgebra} or just a {\it Heisenberg superalgebra}.

For the remainder of the paper, let $\mathfrak{h}$ be  finite-dimensional vector space over $\mathbb{C}$ equipped with a nondegenerate symmetric bilinear form $\langle \cdot , \cdot \rangle$.   Let $d$ denote the dimension of $\mathfrak{h}$, let $t$ and $x$ denote formal commuting variables, and let $U(\cdot)$ denote the universal enveloping algebra for a Lie superalgebra $(\cdot)$.

\subsection{Bosonic Heisenberg VOAs}\label{bosonic-section}

Form the affine Lie algebra
\[\hat{\mathfrak{h}} = \mathfrak{h} \otimes \mathbb{C}[t,t^{-1}] \oplus \mathbb{C} \mathbf{k} ,\]
with the Lie bracket relations 
\begin{eqnarray}
\bigl[\mathbf{k}, \hat{\mathfrak{h}} \bigr] &=& 0 \label{commutation1}\\
\bigl[\alpha \otimes t^m, \beta \otimes t^n \bigr] &=& \langle \alpha, \beta \rangle m \delta_{m + n,0} \mathbf{k} \label{commutation2}
\end{eqnarray}
for $\alpha, \beta \in \mathfrak{h}$ and $m, n \in \mathbb{Z}$.  Then $\hat{\mathfrak{h}}$ is a $\mathbb{Z}$-graded Lie algebra 
\[\hat{ \mathfrak{h} } = \coprod_{n \in \mathbb{Z}}\hat{ \mathfrak{h}}_n \]
where $\hat{\mathfrak{h}}_0 = \mathfrak{h} \oplus \mathbb{C} \mathbf{k}$, and  $\hat{\mathfrak{h}}_n = \mathfrak{h} \otimes t^{-n}$, for $n \neq 0$,
and has graded subalgebras 
\[\hat{\mathfrak{h}}_+ = \mathfrak{h} \otimes t^{-1} \mathbb{C}[t^{-1}] \quad \mathrm{and} \quad \hat{\mathfrak{h}}_- = \mathfrak{h} \otimes t \mathbb{C}[t] .\]
Set 
\[\hat{\mathfrak{h}}_* = \hat{\mathfrak{h}}_- \oplus \hat{\mathfrak{h}}_+ \oplus \mathbb{C} \mathbf{k} = \coprod_{n \neq 0} (\mathfrak{h} \otimes t^n) \oplus \mathbb{C} \mathbf{k} .\]
Then $\hat{\mathfrak{h}}_*$ is a Heisenberg algebra.  In addition $\hat{\mathfrak{h}}_*$ and $ \mathfrak{h}$ are ideals of $\hat{\mathfrak{h}}$ and $\hat{\mathfrak{h}} = \hat{\mathfrak{h}}_* \oplus \mathfrak{h}$.

Let $\mathbb{C}$ be the $(\hat{\mathfrak{h}}_- \oplus \hat{\mathfrak{h}}_0)$-module such that $\hat{\mathfrak{h}}_-$ and $\mathfrak{h}$ act trivially and $\mathbf{k}$ acts as $1$.  Let 
\[V_{bos} =  U(\hat{\mathfrak{h}}) \otimes_{U(\hat{\mathfrak{h}}_- \oplus
\hat{\mathfrak{h}_0}) } \mathbb{C}  \cong S(\hat{\mathfrak{h}}_+) \]
so that $V_{bos}$ is naturally isomorphic to the symmetric algebra of polynomials in $\hat{\mathfrak{h}}_+$; see Remark \ref{symmetric-algebra-remark}.  It is also the universal enveloping algebra for $\hat{\mathfrak{h}}_+$.  Let $\alpha \in \mathfrak{h}$ and $n \in \mathbb{Z}$.  We will use the notation
\[\alpha(n) = \alpha \otimes t^n \in \hat{\mathfrak{h}}.\]

Note that $V_{bos}$ is a $\hat{\mathfrak{h}}$-module with action induced from the commutation relations (\ref{commutation1}) and (\ref{commutation2}) given by
\begin{eqnarray}
\mathbf{k} \beta (-m)\mathbf{1} &=& \beta (-m) \mathbf{1}\label{boson-action1}\\
\alpha(0)  \beta (-m)\mathbf{1} &=& 0\\
\alpha (n) \beta (-m)\mathbf{1} &=& \langle \alpha, \beta \rangle n \delta_{m,n}\mathbf{1} \\
\alpha (-n) \beta (-m)\mathbf{1} &=& \beta(-m) \alpha(-n) \mathbf{1} \label{boson-action-last}
\end{eqnarray}
for $\alpha, \beta \in \mathfrak{h}$ and $m, n \in \mathbb{Z}_+$.

\begin{rema}\label{symmetric-algebra-remark}
{\em
Let  $\{\alpha^{(1)}, \alpha^{(2)}, \dots, \alpha^{(d)} \}$ be an orthonormal basis for $\mathfrak{h}$.  Let $a^{(j)}_n$, for $n \in \mathbb{Z}_+$, be mutually commuting independent formal variables.  Then $\hat{\mathfrak{h}}$ acts on the space 
\begin{equation}\label{symmetric-space}
  \mathbb{C} [a^{(1)}_1, a^{(1)}_2, \dots, a^{(2)}_1, a^{(2)}_2, \dots, a^{(d)}_1, a^{(d)}_2, \dots ]
\end{equation}
by 
\begin{eqnarray}
\mathbf{k} &\mapsto  &  1\\
\alpha(0) &\mapsto& 0\\
\alpha^{(j)} (n)  &\mapsto& n \frac{\partial}{\partial  a^{(j)}_n}\\
\alpha^{(j)} (-n) &\mapsto& a^{(j)}_n, \label{multiplication}
\end{eqnarray}
for $n \in \mathbb{Z}_+$, and where the operator on the left of (\ref{multiplication}) is the multiplication operator.  Then the symmetric algebra (\ref{symmetric-space}) is isomorphic to $V_{bos}$ as an $\hat{\mathfrak{h}}$-module.}
\end{rema} 

For $\alpha \in \mathfrak{h}$, set
\[ \alpha(x)^b = \sum_{n \in \mathbb{Z}} \alpha(n) x^{-n-1} ,\]
Define the {\it normal ordering} operator ${}^\circ_\circ \cdot {}^\circ_\circ$ on products of the operators $\alpha(n)$ by
\[ {}^\circ_\circ\alpha(m) \beta(n) {}^\circ_\circ = \left\{ \begin{array}{ll} \alpha (m) \beta (n) & \mbox{if $m\leq n$}\\
\beta (n) \alpha (m) & \mbox{if $m> n$}\\
\end{array} \right. \]
for $m,n \in \mathbb{Z}$. 

For $n \in \mathbb{N}$, let  
\begin{equation}
\partial_n = \frac{1}{n!} \left(\frac{d}{dx}\right)^n.
\end{equation}
For $v = \alpha_1(-n_1) \alpha_2 (-n_2) \cdots \alpha_m (-n_m) \mathbf{1} \in V_{bos}$, for $\alpha_j \in \mathfrak{h}$, $n_j \in \mathbb{Z}_+$, and $j = 1, \dots, m$ and $m \in \mathbb{N}$, define the vertex operator corresponding to $v$ to be 
\begin{equation}
Y(v, x) = {}^\circ_\circ \left( \partial_{n_1-1} \alpha_1(x)^b \right) \left( \partial_{n_2-1} \alpha_2(x)^b \right) \cdots \left( \partial_{n_m-1} \alpha_m(x)^b \right) {}^\circ_\circ .
\end{equation}

Note that
\begin{equation}
[\alpha^{(j)} (x_1)^b, \alpha^{(k)}(x_2)^b] = \delta_{j,k} \left( \frac{1}{(x_1 - x_2)^2} - \frac{1}{(-x_2 + x_1)^2} \right)
\end{equation}
implying that the $\alpha^{(j)}(x)^b = Y(\alpha^{(j)} (-1)\mathbf{1}, x)$, for $j = 1,\dots, d$, are mutually local.  Setting 
\[ \omega_{bos} = \frac{1}{2} \sum_{j = 1}^d \alpha^{(j)}(-1) \alpha^{(j)}(-1) \mathbf{1},\]
we have that 
\[L(-1) = \sum_{j = 1}^d \sum_{n \in \mathbb{Z}_+}  \alpha^{(j)} (-n) \alpha^{(j)} (n - 1) , \]
and thus
\[ [L(-1), Y(\alpha^{(j)} (-1)\mathbf{1}, x)] =  [L(-1), \alpha^{(j)} (x)^b] = \frac{d}{d x} \alpha^{(j)} (x)^b = \frac{d}{d x } Y(\alpha^{(j)} (-1)\mathbf{1} , x).  \] 
In addition, we have 
\begin{eqnarray}
\omega_{0} \omega &=&  L(-1) \omega \label{Virasoro-condition1} \\
\omega_{1} \omega &=& L(0) \omega  \ = \ 2 \omega \\
\omega_{2} \omega &=&  L(1) \omega  \ \ = \ 0\\
\omega_3 \omega  &=&L(2) \omega \ \ = \   \frac{1}{2} c \mathbf{1} \\
\omega_n \omega &=&  0 \qquad  \mbox{for $n \in \mathbb{Z}$ with $n\geq 4$.} \label{Virasoro-condition-last}
\end{eqnarray}

Since $V_{bos}$ is generated by the $\alpha^{(j)}(-1)$, which we denote by 
\[V_{bos} = \langle \alpha^{(1)}(-1)\mathbf{1}, \dots, \alpha^{(d)} (-1)\mathbf{1}\rangle,\] 
then by for instance \cite{LL}, we have that  $(V_{bos}, Y, \mathbf{1}, \omega_{bos})$ is a vertex operator algebra with central charge  $d$, where of course the vacuum vector is just $1$.  $V_{bos}$ is called the {\it rank $d$ Heisenberg VOA} or the {\it $d$ free boson VOA}.   

Since
\[L(0) =  \sum_{j = 1}^d \sum_{n \in \mathbb{Z}_+} \left(\alpha^{(j)}(-n) \alpha^{(j)}(n) + \frac{1}{2} \alpha^{(j)} (0)^2 \right)\]
the graded dimension of $V_{bos}$ using the $\mathbb{Z}$-grading of $V_{bos}$ by eigenvalues of $L(0)$ is
\begin{equation}
\mathrm{dim}_q V_{bos} = q^{-c/24} \sum_{n \in \mathbb{Z}} \mathrm{dim} (V_{bos})_n q^n = q^{-d/24} \prod_{n \in \mathbb{Z}_+} (1-q^n)^{-d} = \frac{1}{(\eta (q))^d},
\end{equation}
where $\eta(q)$ is the Dedekind $\eta$-function.

\subsection{Free Fermionic VOSAs}\label{fermionic-section}

Form the affine Lie superalgebra
\[\hat{\mathfrak{h}}^f = \mathfrak{h} \otimes t^{1/2}\mathbb{C}[t,t^{-1}]\oplus \mathbb{C} \mathbf{k} ,\]
with $\mathbb{Z}_2$-grading given by $\mathrm{sgn}(\alpha \otimes t^n)= 1$ for $ n \in  \mathbb{Z}+ \frac{1}{2}$, and $\mathrm{sgn}(\mathbf{k})= 0$, and Lie super-bracket relations 
\begin{eqnarray}
\bigl[\mathbf{k}, \hat{\mathfrak{h}}^f \bigr] &=& 0 \label{commutation1F}\\
\bigl[\alpha \otimes t^m, \beta \otimes t^n \bigr] &=& \langle \alpha, \beta \rangle  \delta_{m + n,0} \mathbf{k} \label{commutation2F}
\end{eqnarray}
for $\alpha, \beta \in \mathfrak{h}$ and $m, n \in \mathbb{Z} + \frac{1}{2}$.  Then $\hat{\mathfrak{h}}^f$ is a $(( \mathbb{Z} + \frac{1}{2}) \cup \{0\})$-graded Lie superalgebra 
\[\hat{ \mathfrak{h} }^f = \coprod_{n \in (\mathbb{Z} + \frac{1}{2})\cup \{0\} }\hat{ \mathfrak{h}}^f_n \]
where $\hat{\mathfrak{h}}^f_n = \mathfrak{h} \otimes t^{-n}$, for $n \in \mathbb{Z} + \frac{1}{2}$, and $\hat{\mathfrak{h}}^f_0  = \mathbb{C}\mathbf{k}$. 
It has graded subalgebras 
\[\hat{\mathfrak{h}}^f_+ = \mathfrak{h} \otimes t^{-1/2} \mathbb{C}[t^{-1}] \quad \mathrm{and} \quad \hat{\mathfrak{h}}^f_- = \mathfrak{h} \otimes t^{1/2} \mathbb{C}[t] .\]
Note that 
\[\hat{\mathfrak{h}}^f = \hat{\mathfrak{h}}^f_- \oplus \hat{\mathfrak{h}}^f_+ \oplus \mathbb{C} \mathbf{k},\]
and note that $\hat{\mathfrak{h}}^f$ is a Heisenberg superalgebra. 

Let $\mathbb{C}$ be the $(\hat{\mathfrak{h}}^f_- \oplus \mathbb{C}\mathbf{k})$-module such that $\hat{\mathfrak{h}}^f_-$ acts trivially and $\mathbf{k}$ acts as $1$.  Let 
\[V_{fer} = U(\hat{\mathfrak{h}}^f) \otimes_{U(\hat{\mathfrak{h}}^f_- \oplus \mathbb{C}\mathbf{k})} \mathbb{C}  \cong \mbox{$\bigwedge$}(\hat{\mathfrak{h}}^f_+) ,\]
so that $V_{fer}$ is naturally isomorphic to the algebra of polynomials in the anticommuting elements of $\hat{\mathfrak{h}}^f_+$; see Remark \ref{exterior-algebra-remark}.  

Let $\alpha \in \mathfrak{h}$ and $n \in \mathbb{Z} + \frac{1}{2}$.  We will use the notation
\[\alpha(n) = \alpha \otimes t^n.\]

Then $V_{fer}$ is a $\hat{\mathfrak{h}}^f$-module with action induced from the supercommutation relations (\ref{commutation1F}) and (\ref{commutation2F}) given by
\begin{eqnarray}
\mathbf{k} \beta (-m) \mathbf{1}&=& \beta (-m) \mathbf{1} \label{fermion-action1} \\
\alpha (n) \beta (-m)\mathbf{1} &=& \langle \alpha, \beta \rangle  \delta_{m,n} \mathbf{1} \\
\alpha (-n) \beta (-m)\mathbf{1} &=& - \beta(-m) \alpha(-n)\mathbf{1} \label{fermion-action-last}
\end{eqnarray}
for $\alpha, \beta \in \mathfrak{h}$ and $m, n \in  \mathbb{N} + \frac{1}{2}$.

\begin{rema}\label{exterior-algebra-remark}
{\em
Let  $\{\alpha^{(1)}, \alpha^{(2)}, \dots, \alpha^{(d)} \}$ be an orthonormal basis for $\mathfrak{h}$.  Let $a^{(j)}_n$, for $n \in  \mathbb{N} + \frac{1}{2}$, be formal variables.  Then $\hat{\mathfrak{h}}^f$ acts on the space  
\begin{equation}\label{exterior-space}
\mbox{$\bigwedge$}  \left[a^{(1)}_{\frac{1}{2}}, a^{(1)}_{\frac{3}{2}}, \dots, a^{(2)}_\frac{1}{2}, a^{(2)}_\frac{3}{2}, \dots, a^{(d)}_\frac{1}{2}, a^{(d)}_\frac{3}{2}, \dots \right],
\end{equation}
by
\begin{eqnarray}
\mathbf{k} &\mapsto  &  1\\
\alpha^{(j)} (n)  &\mapsto&  \frac{\partial}{\partial  a^{(j)}_n}\\
\alpha^{(j)} (-n) &\mapsto& a^{(j)}_n, \label{multiplicationF}
\end{eqnarray}
for $n \in \mathbb{N} + \frac{1}{2}$ and where the operator on the left of (\ref{multiplicationF}) is the multiplication operator. That is, if we consider the $a^{(j)}_n$, for $n \in  \mathbb{N} + \frac{1}{2}$, as formal mutually anti-commuting variables,  then the resulting exterior algebra (\ref{exterior-space}) is isomorphic to $V_{fer}$ as an $\hat{\mathfrak{h}}^f$-module.  
}
\end{rema} 

For $\alpha \in \mathfrak{h}$, set
\begin{equation} \alpha(x)^f = \sum_{n \in \frac{1}{2} + \mathbb{Z}} \alpha(n) x^{-n-\frac{1}{2}} ,
\end{equation}
Define the {\it normal ordering} operator ${}^\circ_\circ \cdot {}^\circ_\circ$ on products of the operators $\alpha(n)$ by
\begin{equation}
{}^\circ_\circ \alpha(m) \beta(n) {}^\circ_\circ = \left\{ \begin{array}{ll} \alpha (m) \beta (n) & \mbox{if $m\leq n$}\\
- \beta (n) \alpha (m) & \mbox{if $m> n$}\\
\end{array} \right. 
\end{equation}
for $m,n \in  \mathbb{Z} + \frac{1}{2}$. 

For $v = \alpha_1(-n_1) \alpha_2 (-n_2) \cdots \alpha_m (-n_m) \mathbf{1} \in V_{fer}$, for $\alpha_j \in \mathfrak{h}$, $n_j \in \mathbb{N}+ \frac{1}{2}$, and $j = 1, \dots, m$ and $m \in \mathbb{N}$, define the vertex operator corresponding to $v$ to be 
\begin{equation}
Y(v, x) = {}^\circ_\circ \left( \partial_{n_1-\frac{1}{2}} \alpha_1(x)^f \right) \left( \partial_{n_2-\frac{1}{2}} \alpha_2(x)^f \right) \cdots \left( \partial_{n_m-\frac{1}{2}} \alpha_m(x)^f\right) {}^\circ_\circ .
\end{equation}

Note that
\begin{equation}
[\alpha^{(j)} (x_1)^f, \alpha^{(k)}(x_2)^f] = \delta_{j,k} \left( \frac{1}{(x_1 - x_2)} - \frac{1}{(-x_2 + x_1)} \right)
\end{equation}
implying that the $\alpha^{(j)}(x)^f = Y(\alpha^{(j)} (-1/2)\mathbf{1}, x)$, for $j = 1,\dots, d$, are mutually local.  Setting 
\begin{equation}\label{fermionic-Virasoro}
 \omega_{fer} = \frac{1}{2} \sum_{j = 1}^d \alpha^{(j)}(-3/2) \alpha^{(j)}(-1/2) \mathbf{1},
\end{equation}
we have that 
\[L(-1) = \sum_{j=1}^d \sum_{n \in \mathbb{N}+ \frac{1}{2}} (n - 1/2)  \alpha^{(j)} (-n) \alpha^{(j)} (n-1) ,\]
and thus
\begin{eqnarray*}
[ L(-1) , Y(\alpha^{(j)} (-1/2)\mathbf{1}, x)] &=& [L(-1), \alpha^{(j)} (x)^f] \ = \ \frac{d}{d x} \alpha^{(j)} (x)^f \\
&=& \frac{d}{d x} Y(\alpha^{(j)}(-1/2)\mathbf{1}, x).
\end{eqnarray*}
In addition, we have that the conditions (\ref{Virasoro-condition1})--(\ref{Virasoro-condition-last}) are satisfied for $\omega_{fer}$.
Since $V_{fer} = \langle \alpha^{(1)}(-1/2)\mathbf{1}, \dots, \alpha^{(d)} (-1/2)\mathbf{1}\rangle$, then by for instance \cite{Li-untwisted}, it follows that $(V_{fer}, Y, \mathbf{1}, \omega_{fer})$ is a vertex operator superalgebra with  central charge $d/2$.  $V_{fer}$ is called the {\it $d$ free fermion VOSA}.

When $d$ is even, $V_{fer}$ is precisely the VOSA studied in \cite{FFR} denoted $CM(\mathbb{Z} + \frac{1}{2})$, although in \cite{FFR} a polarized basis for $\mathfrak{h}$ is used;  see Remark \ref{polarization-remark} below.

Since 
\begin{equation}
L(0)_{fer}=  \sum_{j = 1}^d \sum_{n \in \frac{1}{2} + \mathbb{N}}n \alpha^{(j)}(-n) \alpha^{(j)}(n) ,
\end{equation}
the graded dimension of $V_{fer}$ using the $\mathbb{Z}$-grading of $V_{fer}$ by eigenvalues of $L(0)_{fer}$ is
\begin{equation}
\mathrm{dim}_q V_{fer} = q^{-c/24} \sum_{n \in \frac{1}{2}\mathbb{Z}} \mathrm{dim} (V_{fer})_n q^n = q^{-d/48} \prod_{n \in \mathbb{Z}_+} (1+q^{n-1/2})^{d} = \mathfrak{f}(q)^d,
\end{equation}
where $\mathfrak{f}(q)$ is a classical Weber function \cite{YZ}.   A simple calculation shows that in fact $\mathfrak{f}(q) = \frac{\eta(q)^2}{\eta(q^2) \eta(q^{1/2})}$.   

In addition, the {\it superdimension} of a vertex operator superalgebra $V = V^0 \oplus V^1$ is sometimes of interest.  It is defined to be 
\begin{equation}
\mathrm{sdim}_q V = \mathrm{dim}_q V^{(0)} - \mathrm{dim}_q V^{(1)}.
\end{equation}
Thus the superdimension of $V_{fer}$ is
\begin{equation}
\mathrm{sdim}_q V_{fer} = q^{-d/48} \prod_{n \in \mathbb{Z}_+} (1-q^{n-1/2})^{d}  = \mathfrak{f}_1(q)^d
\end{equation}
where $\mathfrak{f}_1(q)$ is also a classical Weber function.  Observe that $\mathfrak{f}_1(q) = \frac{\eta(q^{1/2})}{\eta(q)}$.  

\begin{rema}\label{Weber-remark}{\em
In addition to the two classical Weber functions, $\mathfrak{f}$ and $\mathfrak{f}_1$, there is a third classical Weber function, denoted $\mathfrak{f}_2$ and given by 
\begin{equation}
\mathfrak{f}_2(q) = \sqrt{2} q^{1/24} \prod_{n \in \mathbb{Z}_+} (1+ q^n) = \sqrt{2} \frac{\eta(q^2)}{\eta(q)}.
\end{equation}
This third classical Weber function, $\mathfrak{f}_2$, will appear in Section \ref{N=1,2-Ramond-section}.  These three Weber functions, $\mathfrak{f}, \mathfrak{f}_1$, and $\mathfrak{f}_2$, form a set that is $SL_2(\mathbb{Z})$-invariant up to permutation and multiplication by 48-th roots of unity \cite{YZ}.
}
\end{rema}

\subsection{Free Boson-Fermion N=1 VOSAs}\label{boson-fermion-section}

Let 
\[V = V_{bos} \otimes V_{fer}\]
and set $Y(u \otimes v, x) = Y(u,x)\otimes Y(v,x)$ for $u \in V_{bos}$ and $v \in V_{fer}$.  Then essentially by \cite{FHL}, $V$ is a VOSA with $\omega = \omega_{bos} \otimes \mathbf{1} + \mathbf{1} \otimes \omega_{fer}$ with central charge $c = 3d/2$ and with graded dimension
\begin{equation}
\mathrm{dim}_q V =  \left(\frac{\mathfrak{f}(q)}{\eta(q)}\right)^d  = \left(\frac{\eta(q)}{\eta(q^2) \eta(q^{1/2})} \right)^d .
\end{equation}
The superdimension is given by
\begin{equation}
\mathrm{sdim}_q V =  \left(\frac{\mathfrak{f}_1 (q)}{\eta(q)}\right)^d  = \left(\frac{\eta(q^{1/2})}{\eta(q)^2} \right)^d .
\end{equation}

Since  $\omega = \omega_{bos} \otimes \mathbf{1} + \mathbf{1} \otimes \omega_{fer}$, we have that 
\begin{equation}
L(n) = \frac{1}{2} \sum_{j = 1}^d  \sum_{m  \in \mathbb{Z}} \left( {}^\circ_\circ \alpha^{(j)}( m) \alpha^{(j)} (n- m ) -m   \alpha^{(j)}( m - 1/2) \alpha^{(j)} (n- m + 1/2)  {}^\circ_\circ \right),
\end{equation}
where we are suppressing the tensor product symbol.  
In addition, setting 
\begin{equation}\label{NS-element}
\tau = \sum_{j =1}^d \alpha^{(j)}(-1) \alpha^{(j)}(-1/2) \mathbf{1},
\end{equation}
we have that
\begin{equation}
\tau_{n+1} = G(n + 1/2) = \sum_{j = 1}^d  \sum_{m \in  \mathbb{Z} + \frac{1}{2}} \alpha^{(j)}( m) \alpha^{(j)} (n- m + 1/2)
\end{equation}
and
\begin{eqnarray}
\tau_0 \tau &=& G(-1/2) \tau \ = \ 2 \omega \\
\tau_1 \tau &=& G(1/2) \tau \ = \ 0 \\
\tau_2 \tau &=&  G(3/2) \tau \ = \ \frac{2}{3} c \mathbf{1} \\
\tau_n \tau &=& 0 \qquad  \mbox{for $n \in \mathbb{Z}$ with $n\geq 3$}
\end{eqnarray}
hold.  This implies that $\omega_{n} = L(n -1)$ (which satisfy (\ref{Virasoro-condition1})--(\ref{Virasoro-condition-last})) and $\tau_{n + 1} = G(n + 1/2)$, for $n \in \mathbb{Z}$, satisfy the N=1 Neveu-Schwarz relations with central charge $c \in \mathbb{C}$, which in this case is $3d/2$.

It follows that $(V, Y, \mathbf{1}, \tau)$ is an N=1 VOSA.  We will call $V$ the {\it $d$ free boson-fermion N=1 VOSA}.

In particular, 
\begin{equation}
\tau_0 = G(-1/2) = \sum_{j = 1}^d  \sum_{m \in \frac{1}{2} \mathbb{Z}_+} \alpha^{(j)}( -m) \alpha^{(j)} (m-1/2)
\end{equation}
and letting $\varphi$ be an odd formal variable, then setting
\begin{equation}
Y(u, (x, \varphi)) = Y(u,x) + \varphi Y(G(-1/2)u, (x, \varphi)) 
\end{equation}
we recover the odd component of the vertex operators; see \cite{B-announce}, \cite{B-vosas}.

Thus setting
\[\alpha(x, \varphi) = \alpha(x)^f + \varphi \alpha(x)^b ,\]
for $\alpha \in \mathfrak{h}$, we have for  $v = \alpha_1(-n_1) \alpha_2 (-n_2) \cdots \alpha_k (-n_k) \mathbf{1} \in V$, for $\alpha_j \in \mathfrak{h}$, $n_j \in \frac{1}{2} \mathbb{Z}_+$, and $j = 1, \dots, k$ and $k \in \mathbb{N}$, the vertex operator corresponding to $v$ is given by
\begin{multline}
Y(v, (x, \varphi)) = {}^\circ_\circ \left( \frac{1}{(\lfloor n_1-\frac{1}{2}\rfloor) !} D^{2n_1-1} \alpha_1(x, \varphi)  \right) \\
\left( \frac{1}{(\lfloor n_2-\frac{1}{2} \rfloor)!}D^{2n_2-1} \alpha_2(x, \varphi) \right) 
 \cdots \left( \frac{1}{(\lfloor n_k-\frac{1}{2} \rfloor)!} D^{2n_k-1} \alpha_k(x, \varphi)\right) {}^\circ_\circ
\end{multline}
where $\lfloor \cdot \rfloor$ is the floor function and $D = \frac{\partial}{\partial \varphi} + \varphi \frac{\partial}{\partial x}$ is the odd superderivation satisfying $D^2 = \frac{\partial}{\partial x}$.  

\subsection{Free  N=2 VOSAs}\label{N=2-free-section}

Again form the free boson-fermion N=1 VOSA, $V= V_{bos} \otimes V_{fer} =  S( \hat{\mathfrak{h}}_+) \otimes \bigwedge (\hat{\mathfrak{h}}_+^f)$  as in Section \ref{boson-fermion-section}.  Now take the tensor product of two copies of $V$, i.e,  $V \otimes V$. 

As in Section \ref{boson-fermion-section}, letting $\{\alpha^{(j)} \;| \; j = 1, \dots, d \}$ be an orthonormal basis for $\mathfrak{h}$, then  $V\otimes V$ is an N=1 VOSA with N=1 superconformal element 
\begin{equation}
\tau^{(1)} = \sum_{j=1}^d \left( \alpha^{(j)}(-1) \alpha^{(j)}(-1/2)\mathbf{1} \otimes \mathbf{1} + \mathbf{1} \otimes  \alpha^{(j)}(-1) \alpha^{(j)}(-1/2)\mathbf{1} \right),
\end{equation}
central charge $c = 3d$, and conformal element
\begin{multline}
\omega = \frac{1}{2} \sum_{j=1}^d \left( \alpha^{(j)}(-1) \alpha^{(j)}(-1)\mathbf{1} \otimes \mathbf{1} + \mathbf{1} \otimes \alpha^{(j)}(-1) \alpha^{(j)}(-1)\mathbf{1} \right. \\
 \left. +   \alpha^{(j)}(-3/2) \alpha^{(j)}(-1/2)\mathbf{1} \otimes \mathbf{1} + \mathbf{1} \otimes  \alpha^{(j)}(-3/2) \alpha^{(j)}(-1/2)\mathbf{1}\right).
\end{multline}
But in addition, $V \otimes V$ is an N=2 VOSA with the weight-one vector $\mu$ giving the $J(n)$, for $n\in\Z$, (which generate a representation of affine $\mathfrak{u}(1)$) given by
\begin{equation}
\mu = i \sum_{j=1}^d \alpha^{(j)}(-1/2)\mathbf{1} \otimes \alpha^{(j)}(-1/2)\mathbf{1} ,
\end{equation}
and the other N=1 superconformal element given by
\begin{equation}
\tau^{(2)} = \sum_{j=1}^d \left( \alpha^{(j)}(-1)\mathbf{1} \otimes  \alpha^{(j)}(-1/2) \mathbf{1} -  \alpha^{(j)}(-1/2)\mathbf{1} \otimes \alpha^{(j)}(-1)\mathbf{1} \right).
\end{equation}

For $\alpha \in \mathfrak{h}$, we write $\alpha_{(1)} (n) = \alpha(n) \otimes 1$ and $\alpha_{(2)} (n) = 1 \otimes \alpha(n)$ for $n \in \frac{1}{2} \mathbb{Z}$.  Therefore setting
\begin{eqnarray}
\qquad \alpha_{(1)}(x, \varphi_1, \varphi_2) \! \! \! &=& \! \! \!  \alpha_{(1)}(x)^f + \varphi_1 \alpha_{(1)}(x)^b - \varphi_2  \alpha_{(2)}(x)^b + \varphi_1 \varphi_2 \frac{\partial}{\partial x} \alpha_{(2)}(x)^f ,\\
\alpha_{(2)}(x, \varphi_1, \varphi_2) \! \! \! &=& \! \! \! \alpha_{(2)}(x)^f + \varphi_1 \alpha_{(2)}(x)^b + \varphi_2 \alpha_{(1)}(x)^b - \varphi_1 \varphi_2 \frac{\partial}{\partial x} \alpha_{(1)}(x)^f,
\end{eqnarray}
for $\alpha \in \mathfrak{h}$, then for 
$v = \alpha_{1, (k_1)} (-n_1) \alpha_{2, (k_2)} (-n_2) \cdots \alpha_{l, (k_l)} (-n_l) \cdot ( \mathbf{1} \otimes \mathbf{1}) \in V\otimes V$, for $\alpha_{m} \in \mathfrak{h}$, $n_m \in \frac{1}{2} \mathbb{Z}_+$, $m = 1, \dots, l$, $k_m = 1,2$, and $l \in \mathbb{N}$, the vertex operator corresponding to $v$ is given by
\begin{eqnarray}
\qquad Y(v, (x, \varphi_1, \varphi_2)) \! \! \! &=& \! \! \!  {}^\circ_\circ \left( \frac{1}{(\lfloor n_1-\frac{1}{2}\rfloor) !} (D^{(k_1)})^{2n_1-1} \alpha_{1, (k_1)}(x, \varphi_1, \varphi_2)  \right)\\
& & \cdot  \left( \frac{1}{(\lfloor n_2-\frac{1}{2} \rfloor)!}(D^{(k_2)})^{2n_2-1} \alpha_{2,(k_2)}(x, \varphi_1, \varphi_2) \right) \nonumber \\
& &  \cdots \left( \frac{1}{(\lfloor n_m-\frac{1}{2} \rfloor)!}( D^{(k_m)})^{2n_m-1} \alpha_{m, (k_m)}(x, \varphi_1, \varphi_2)\right)  {}^\circ_\circ \nonumber
\end{eqnarray}
where $D^{(k_j)} = \frac{\partial}{\partial \varphi_{k_j}} + \varphi_{k_j} \frac{\partial}{\partial x}$ and $\lfloor \cdot \rfloor$ is the floor function. 

Transforming to the homogeneous basis, we set
\begin{equation}
\alpha^{\pm} = \frac{1}{\sqrt{2}} (\alpha_{(1)} \mp i \alpha_{(2)})
\end{equation} 
or equivalently
\begin{equation}
\alpha_{(1 )}= \frac{1}{\sqrt{2}} (\alpha^+ + \alpha^-) \qquad \mbox{and} \qquad \alpha_{(2)} = \frac{i}{\sqrt{2}} (\alpha^+ - \alpha^-).
\end{equation}
Then conformal and superconformal elements in the homogeneous basis are given by
\begin{eqnarray}
\mu &=& \sum_{j = 1}^d \alpha^{(j), +}(-1/2) \alpha^{(j), -}(-1/2) \mathbf{1} \\
\qquad \tau^{(\pm)} &=& \sqrt{2} \sum_{j = 1}^d \alpha^{(j), \mp}(-1) \alpha^{(j), \pm}(-1/2) \mathbf{1} \\
\omega &=& \frac{1}{2} \sum_{j = 1}^d \left( 2 \alpha^{(j), +}(-1) \alpha^{(j), -}(-1) \mathbf{1}+ \alpha^{(j), +}(-3/2) \alpha^{(j), -}(-1/2) \mathbf{1} \right. \\
& & \left. \quad + \, \alpha^{(j), -}(-3/2) \alpha^{(j), +}(-1/2) \mathbf{1}\right). \nonumber
\end{eqnarray}

\begin{rema}{\em
From Remarks \ref{symmetric-algebra-remark} and \ref{exterior-algebra-remark}, we have that as an $(\hat{\mathfrak{h}} \otimes \hat{\mathfrak{h}}^f)^{\otimes 2}$-module, 
\begin{equation}\label{sym-ext-space}
V \otimes V\cong  \mathbb{C} \left[a^{(j), +}_n, a^{(j), -}_n \, | \, j = 1,\dots, d, \ n \in \frac{1}{2}\mathbb{Z}_+ \right]
\end{equation}
where the $a^{(j), \pm}_n$, for $n \in \frac{1}{2} \mathbb{Z}_+$, are commuting formal variables if $n \in \mathbb{Z}_+$ and anti-commuting if $n \in \mathbb{N} + \frac{1}{2}$, and we have the following operators on (\ref{sym-ext-space})
\begin{eqnarray}
\mathbf{k} &\mapsto  &  1\\
\alpha^{(j), \pm} (n)  &\mapsto& n \frac{\partial}{\partial  a^{(j), \mp}_n} \qquad \mbox{for $n \in \mathbb{Z}_+$}\\
\alpha^{(j), \pm} (n)  &\mapsto& \frac{\partial}{\partial  a^{(j), \mp}_n} \qquad \mbox{for $n \in \mathbb{N} +  \frac{1}{2}$}\\
\alpha^{(j), \pm} (-n) &\mapsto& a^{(j), \pm}_n \qquad \mbox{for $n \in \frac{1}{2}  \mathbb{Z}_+$}, \label{multiplication-homo}
\end{eqnarray}
and where the operator on the left of (\ref{multiplication-homo}) is the multiplication operator.}
\end{rema} 

Observe that the $q$-dimension for the N=2 VOSA $V\otimes V$ is just the square of the $q$-dimension for $V$.   But we also have the $p,q$-dimension, i.e., the dimension graded in terms of eigenvalues of both the $L(0)$ and the $J(0)$ operators, and this is given by
\begin{equation}\label{pq-dimension}
\mathrm{dim}_{p,q}  V \otimes V = q^{-d/24} \eta(q)^{-2d} \prod_{n \in \mathbb{Z}_+} (1 + p q^{n - 1/2})^d (1 + p^{-1} q^{n - 1/2})^d.
\end{equation}
Note that equation (\ref{pq-dimension}) contains the Jacobi Triple Product Identity, cf. \cite{Berndt}.

\subsection{Extensions to lattice N=1 and N=2 VOSAs}\label{lattice-section}

Suppose $L$ is a positive definite integral lattice.  Then letting $\mathfrak{h} = L \otimes_{\mathbb{Z}} \mathbb{C}$, and following, for instance \cite{LL}, \cite{S}, \cite{Xu}, let $V_L$ be the lattice VOSA.  If $L$ is a positive definite even lattice, then $V_L$ is a vertex operator algebra.

Consider $V_L \otimes V_{fer}$, where $V_{fer}$ is the free fermionic VOSA based on $\mathfrak{h} = L \otimes_{\mathbb{Z}} \mathbb{C}$, and thus $d = \mathrm{rank} \, L$.  Then $V_L \otimes V_{fer}$ is an N=1 VOSA with N=1 Neveu-Schwarz element given by (\ref{NS-element}).  Similarly, tensoring this N=1 VOSA, $V_L \otimes V_{fer}$, with itself, we get an N=2 VOSA.   The $q$-dimension of the lattice N=1 VOSA $V_L \otimes V_{fer}$ is the $q$-dimension of $V_{bos} \otimes V_{fer}$ multiplied by the theta-function for the lattice, $\Theta(L)$.  If $L$ is even, the $p,q$-dimension of the lattice N=2 VOSA $(V_L \otimes V_{fer})^{\otimes 2}$ is the $p,q$-dimension of $(V_{bos} \otimes V_{fer})^{\otimes 2}$ multiplied by $\Theta(L)^2$.  If $L$ is integral, then the $J(0)$-grading is derived from $J(0) e^{\alpha^\pm} = \pm e^{\alpha^\pm}$ for $\alpha \in L$.

\section{Ramond twisted sectors for free and lattice N=1 and N=2 VOSAs}\label{ramond-construction-section}

Following \cite{FFR}, \cite{Li-twisted}, \cite{DZ}, we present the $\sigma$-twisted modules for the free fermion VOSAs constructed in Section \ref{fermionic-section}.  Then we show how these give $\sigma$-twisted modules for the free N=1 and N=2 VOSAs constructed in Sections \ref{boson-fermion-section} and \ref{N=2-free-section}.   We point out that this construction easily extends to all N=1 and N=2 VOSAs for the positive definite lattice cases.  If the lattice is positive definite integral, then we discuss what constructions can be achieved and which are still as yet not known.

\subsection{$\sigma$-twisted sectors for free fermion VOSAs}\label{Ramond-section}

Form the affine Lie superalgebra
\[\hat{\mathfrak{h}}^f[\sigma] = \mathfrak{h} \otimes \mathbb{C}[t,t^{-1}]\oplus \mathbb{C} \mathbf{k} ,\]
with $\mathbb{Z}_2$-grading given by $\mathrm{sgn}(\alpha \otimes t^n)= 1$ for $n \in \mathbb{Z}$, and $\mathrm{sgn}(\mathbf{k})= 0$, and Lie super-bracket relations 
\begin{eqnarray}
\bigl[\mathbf{k}, \hat{\mathfrak{h}}^f[\sigma] \bigr] &=& 0 \label{commutation1F-twisted}\\
\bigl[\alpha \otimes t^m, \beta \otimes t^n \bigr] &=& \langle \alpha \otimes t^m, \beta \otimes t^n \rangle \ = \ \langle \alpha, \beta \rangle  \delta_{m + n,0} \mathbf{k} \label{commutation2F-twisted}
\end{eqnarray}
for $\alpha, \beta \in \mathfrak{h}$ and $m, n \in \mathbb{Z}$, where we have extended the nondegenerate symmetric bilinear form on $\mathfrak{h}$ to a symmetric nondegenerate bilinear form on $\mathfrak{h} \otimes \mathbb{C}[t, t^{-1}]$.

Then $\hat{\mathfrak{h}}^f[\sigma]$ is a $\mathbb{Z}$-graded Lie superalgebra 
\[\hat{ \mathfrak{h} }^f[\sigma] = \coprod_{n \in\mathbb{Z}}\hat{ \mathfrak{h}}^f[\sigma]_n \]
where $\hat{\mathfrak{h}}^f[\sigma]_0 = \mathfrak{h} \oplus \mathbb{C} \mathbf{k}$, and  $\hat{\mathfrak{h}}^f[\sigma]_n = \mathfrak{h} \otimes t^{-n}$ for $n \neq 0$.  And $\hat{\mathfrak{h}}^f[\sigma]$ is a Heisenberg superalgebra.

If $\mathrm{dim} \, \mathfrak{h} = d$ is even, i.e. $d = 2l$, then we can choose a polarization of $\mathfrak{h}$ into maximal isotropic subspaces $\mathfrak{a}^\pm$.  That is $\mathfrak{a}^\pm$ both have dimension $l$, and satisfy $\langle \mathfrak{a}^+, \mathfrak{a}^+ \rangle = \langle \mathfrak{a}^-, \mathfrak{a}^-\rangle = 0$, and we can choose a basis of $\mathfrak{a}^-$, given by $\{\beta^{(1)}_-, \beta^{(2)}_-, \dots, \beta^{(l)}_-\}$, and a dual basis for $\mathfrak{a}^+$, given by $\{\beta^{(1)}_+, \beta^{(2)}_+, \dots, \beta^{(l)}_+\}$ such that $\langle \beta^{(j)}_-, \beta^{(n)}_+ \rangle = \delta_{j,n}$.

 If $\mathrm{dim} \, \mathfrak{h} = d$ is odd, i.e. $d = 2l+1$, then we can choose a polarization of $\mathfrak{h}$ into maximal isotropic subspaces $\mathfrak{a}^\pm$, each of  dimension $l$, and a one-dimensional space $\mathfrak{e}$, so that $\mathfrak{h} = \mathfrak{a}^- \oplus \mathfrak{a}^+ \oplus \mathfrak{e}$, and such that $\langle \mathfrak{a}^\pm, \mathfrak{e}\rangle = 0$, and $\mathfrak{e} = \mathbb{C} \epsilon$ with $\langle \epsilon , \epsilon \rangle = 2$. 

\begin{rema}\label{polarization-remark} \em{If $\{ \alpha^{(1)}, \alpha^{(2)}, \dots, \alpha^{(d)} \}$ is an orthonormal basis for $\mathfrak{h}$ with respect to the symmetric bilinear form as in Section \ref{free-section}, in particular see Remarks \ref{symmetric-algebra-remark} and \ref{exterior-algebra-remark}, then a polarization for $\mathfrak{h}$ can be given as follows: For
$d$ either $2l$ or $2l +1$, set 
\begin{equation}\label{polarize}
\beta^{(j)}_\pm = \frac{1}{\sqrt{2}} \left( \alpha^{(j)} \pm i \alpha^{(j + l)} \right)
\end{equation} 
for $j =1,2,\dots l$.  Then $\mathfrak{a}^\pm = \mathrm{span}_\mathbb{C} \{  \beta^{(1)}_\pm, \beta^{(2)}_\pm, \dots, \beta_\pm^{(l)} \}$ gives a decomposition into maximal polarized spaces.  If $d = 2l+1$, then set $\epsilon = \sqrt{2} \alpha^{(d)}$.   Note that (\ref{polarize}) is equivalent to $\alpha^{(j)} = \frac{1}{\sqrt{2}} \left( \beta_+^{(j)} + \beta_-^{(j)}\right)$ and  $\alpha^{(j+l)} = \frac{-i}{\sqrt{2}} \left( \beta_+^{(j)} - \beta_-^{(j)}\right)$ for $j = 1,\dots,l$.
}
\end{rema}

Then $\hat{\mathfrak{h}}^f[\sigma]$ has the following graded subalgebras
\[\hat{\mathfrak{h}}^f [\sigma]_+ = \mathfrak{h} \otimes t^{-1} \mathbb{C}[t^{-1}] \qquad \mbox{and} \qquad \hat{\mathfrak{h}}^f [\sigma]_- = \mathfrak{h} \otimes t \mathbb{C}[t], \]
and we have $\hat{\mathfrak{h}}^f[\sigma] = \hat{\mathfrak{h}}^f[\sigma]_- \oplus \mathfrak{h} \oplus \hat{\mathfrak{h}}^f[\sigma]_+ \oplus \mathbb{C}\mathbf{k}$.  In addition, $\hat{\mathfrak{h}}^f[\sigma]$ has the subalgebras
\[\hat{\mathfrak{h}}^f[\sigma]_+ \oplus \mathfrak{a}^+ \qquad \mbox{and} \qquad  \hat{\mathfrak{h}}^f[\sigma]_- \oplus \mathfrak{a}^-\]
for $d$ even and
\[\hat{\mathfrak{h}}^f[\sigma]_+ \oplus \mathfrak{a}^+ \oplus \mathfrak{e} \qquad \mbox{and} \qquad  \hat{\mathfrak{h}}^f[\sigma]_- \oplus \mathfrak{a}^-\]
for $d$ odd.

Let $\mathbb{C}$ be the $(\hat{\mathfrak{h}}^f [\sigma]_- \oplus \mathfrak{a}^- \oplus  \mathbb{C} \mathbf{k})$-module such that $\hat{\mathfrak{h}}^f[\sigma]_-\oplus \mathfrak{a}^-$ acts trivially and $\mathbf{k}$ acts as $1$.  Set
\begin{equation}
M_\sigma =  U(\hat{\mathfrak{h}}^f[\sigma]) \otimes_{U(\hat{\mathfrak{h}}^f [\sigma]_- \oplus  \mathfrak{a}^- \oplus \mathbb{C} \mathbf{k})} \mathbb{C}.
\end{equation}
Then as a vector space, we have
\begin{equation}
M_\sigma \stackrel{\mathrm{vec.sp.}}{\cong}\left\{ \begin{array}{ll}
\mbox{$\bigwedge$}(\hat{\mathfrak{h}}^f [\sigma]_+ \oplus \mathfrak{a}^+) & \mbox{if $d$ is even}\\
\mbox{$\bigwedge$}(\hat{\mathfrak{h}}^f [\sigma]_+ \oplus \mathfrak{a}^+ \oplus  \mathfrak{e}) & \mbox{if $d$ is odd}
\end{array} \right.,
\end{equation}
where if $d$ is even, this is also an associative algebra isomorphism, but if $d$ is odd it is not;  rather, if $d$ is odd, $M_\sigma$ is a Clifford algebra but not an exterior algebra.  See Remark \ref{twisted-algebra-remark}.

Let $\alpha \in \mathfrak{h}$ and $n \in \mathbb{Z}$.  We use the notation 
\[ \overline{\alpha(n)} = \alpha \otimes t^n \in \hat{\mathfrak{h}}^f[\sigma] \]
where the overline is meant to distinguish elements of $\hat{\mathfrak{h}}^f[\sigma]$ from elements of $\hat{\mathfrak{h}}$, used to construct the free bosonic theory.

Then $M_\sigma$ is a $\hat{\mathfrak{h}}^f[\sigma]$-module.  For $d$ even, the action induced from the supercommutation relations (\ref{commutation1F-twisted}) and (\ref{commutation2F-twisted}) is given by 
\begin{eqnarray}
\mathbf{k} \overline{\beta (-m)} \mathbf{1} &=& \overline{\beta (-m)} \mathbf{1} \label{twisted-relations1}\\
\overline{\alpha  (n)} \, \overline{\beta (-m) } \mathbf{1} &=& \langle \alpha, \beta \rangle  \delta_{m,n} \mathbf{1} \\
\overline{\alpha (-n)} \, \overline{\beta (-m)} \mathbf{1} &=& - \overline{\beta(-m)} \, \overline{\alpha(-n)}  \mathbf{1}
\end{eqnarray}
for either (i) $\alpha, \beta \in \mathfrak{h}$ and $m, n \in \mathbb{Z}_+$; (ii) $\alpha \in \mathfrak{h}$, $\beta \in \mathfrak{a}^+$, $m =0$, and $n \in \mathbb{Z}_+$;  or (iii) $\alpha \in \mathfrak{a}^-$, $\beta \in \mathfrak{h}$, $n = 0$, and $m \in \mathbb{Z}_+$; and
\begin{equation}\label{twisted-relations-last}
\overline{\alpha(0)} \, \overline{\beta(0)} \mathbf{1} \ = \ \langle \alpha, \beta \rangle \mathbf{1}
\end{equation}
if $\alpha \in \mathfrak{a}^-$ and $\beta \in \mathfrak{a}^+$, and where here $\mathbf{1} = \mathbf{1}_{M_\sigma} = 1$.  

For $d$ odd, the action induced from the supercommutation relations are given by (\ref{twisted-relations1})--(\ref{twisted-relations-last}) as well as  
\begin{eqnarray}
\mathbf{k} \overline{\epsilon(0)} \mathbf{1} &=& \overline{\epsilon(0)}\mathbf{1} \\
\overline{\alpha(0)} \, \overline{\epsilon(0)} \mathbf{1} &=& \langle \alpha, \epsilon \rangle\mathbf{1}
\end{eqnarray}
for $\alpha \in \mathfrak{h}$.

\begin{rema}\label{twisted-algebra-remark}
{\em Let $\{\beta^{(1)}_\pm, \beta^{(2)}_\pm, \dots, \beta^{(l)}_\pm\}$ be the bases for the polarization spaces $\mathfrak{a}^\pm$ as defined in Remark \ref{polarization-remark}.  And, 
if $d$ is odd, let $\mathfrak{e} = \mathbb{C} \epsilon$ with $\langle \epsilon,\epsilon \rangle = 2$.  Then setting $b^{(j)}_{-, n} = \overline{\beta^{(j)}_-(-n)}\mathbf{1}$, for $n \in \mathbb{Z}_+$,  $b^{(j)}_{+,n} = \overline{\beta^{(j)}_+ (n)}\mathbf{1}$ for $n \in \mathbb{N}$, and in addition, if $d$ is odd, $e_n = \overline{\epsilon(-n)} \mathbf{1}$ for $n \in \mathbb{N}$, we have 
\begin{eqnarray*}
M_\sigma &=& \mbox{$\bigwedge$} [\overline{\beta^{(j)}_-(-m)}\mathbf{1}, \, \overline{\beta^{(j)}_+(-n)} \mathbf{1} \, | \, \mbox{for $m \in \mathbb{Z}_+$, $n \in \mathbb{N}$, and $j = 1, \dots, l$} ] \\
&=&  \mbox{$\bigwedge$}  [b^{(1)}_{-,1}, b^{(1)}_{-,2}, \dots, b^{(2)}_{-,1}, b^{(2)}_{-,2}, \dots, b^{(l)}_{-,1}, b^{(l)}_{-,2}, \dots, b^{(1)}_{+,0}, b^{(1)}_{+, 1}, b^{(1)}_{+, 2}, \dots, b^{(2)}_{+,0},  b^{(2)}_{+,1}, \\
& & \quad b^{(2)}_{+,2}, \dots, b^{(l)}_{+,0}, b^{(l)}_{+,1}, b^{(l)}_{+,2}, \dots ],
\end{eqnarray*}
for $d$ even, and in this case, the identification is as an associative algebra.  However, for $d$ odd, we have
\begin{eqnarray*}
M_\sigma &=& \mbox{$\bigwedge$} [\overline{\beta^{(j)}_-(-m)}\mathbf{1}, \, \overline{\beta^{(j)}_+(-n)}\mathbf{1}, \, \overline{\epsilon(-n)}\mathbf{1} \, | \, \mbox{for $m \in \mathbb{Z}_+$, $n \in \mathbb{N}$, and $j = 1, \dots, l$}] \\
&=&  \mbox{$\bigwedge$}  [b^{(1)}_{-,1}, b^{(1)}_{-,2}, \dots, b^{(2)}_{-,1}, b^{(2)}_{-,2}, \dots, b^{(l)}_{-,1}, b^{(l)}_{-,2}, \dots, b^{(1)}_{+,0}, b^{(1)}_{+, 1}, b^{(1)}_{+, 2}, \dots, b^{(2)}_{+,0},  b^{(2)}_{+,1}, \\
& & \quad b^{(2)}_{+,2}, \dots, b^{(l)}_{+,0}, b^{(l)}_{+,1}, b^{(l)}_{+,2}, \dots,
e_0, e_1, e_2, \dots],
\end{eqnarray*}
where in this case, the identification is as a vector space but not as an associative algebra.  As an associative algebra with identity, $M_\sigma$ for $d$ odd is the Clifford algebra generated by $\hat{\mathfrak{h}}^f [\sigma]_+ \oplus \mathfrak{a}^+ \oplus  \mathfrak{e}$ with the corresponding symmetric bilinear form.  

That is, one can think of the $b_{\pm,n}^{(j)}$ as anti-commuting formal variables, and in the case of $d$ even, e.g. $d = 2l$, then $M_\sigma$ is $\mathbb{C}$ adjoin these anti-commuting formal variables.   Then we have the following operators on $M_\sigma$ 
\begin{eqnarray}
\mathbf{k} &\mapsto  \  1  \qquad \ \ \ \ \ \ \ \ \ \ \ \ \ \ \ \ \ \ \ \ \ \ \ \ \ \ &\\
\overline{\beta^{(j)}_\pm (n) } &\mapsto  \ \displaystyle{ \frac{\partial}{\partial  b^{(j)}_{\mp,n}}} \qquad \mathrm{and} \qquad \overline{\beta^{(j)}_- (0) } &\mapsto \ \ \frac{\partial}{\partial  b^{(j)}_{+,0}}\\
\overline{\beta^{(j)}_\pm (-n) }& \! \mapsto \ b^{(j)}_{\pm,n}, \qquad \mathrm{and} \qquad \overline{\beta^{(j)}_+ (0)} &\mapsto \ \ b^{(j)}_{+,0}  \label{multiplication-twisted}
\end{eqnarray}
for $j = 1, \dots, l$, and $n \in \mathbb{Z}_+$, and where the operators on the left of each of the equations in  (\ref{multiplication-twisted}) are multiplication operators.  In addition, if $d$ is odd, e.g. $d = 2l+1$, then $M_\sigma$ also contains the variables $e_{m}$ for $m \in \mathbb{N}$ which anti-commute with the $b_{\pm,n}^{(j)}$, and satisfy
\begin{equation}
e_m e_p = \left\{ \begin{array}{ll}
1 & \mbox{if $m=p=0$}\\
- e_p e_m & \mbox{otherwise}
\end{array} \right. ,
\end{equation}
for $m,p \in \mathbb{N}$.  And we have the following additional operators on $M_\sigma$
\begin{eqnarray}
\overline{\epsilon (n)} &\mapsto & 2\frac{\partial}{\partial e_n}\\
\overline{\epsilon (-m)} &  \mapsto & e_m \label{multiplication-twisted-odd}
\end{eqnarray}
for $n \in \mathbb{Z}_+$, $m \in \mathbb{N}$, and where the operators on the left of  (\ref{multiplication-twisted-odd}) are multiplication operators.
}
\end{rema}

For $\alpha \in \mathfrak{h}$, set
\begin{equation} \alpha(x)^\sigma = \sum_{n \in \mathbb{Z}} \overline{\alpha (n) } x^{-n-\frac{1}{2}} .
\end{equation}
Then for the orthonormal basis of $\mathfrak{h}$, $\alpha^{(j)}$, for $j = 1, \dots, d$, we have 
\begin{equation}
[\alpha^{(j)} (x_1)^\sigma, \alpha^{(k)}(x_2)^\sigma ] = \delta_{j,k} \, x_1^{1/2}x_2^{-1/2}\left( \frac{1}{(x_1 - x_2)} - \frac{1}{(-x_2 + x_1)} \right)\\
\end{equation}
for $j,k = 1,\dots,d$ implying that the $\alpha^{(j)}(x)^\sigma$, for $j = 1,\dots, d$, are mutually local. 

For $v \in V_{fer}$, define $Y^\sigma(v,x) : M_\sigma \longrightarrow M_\sigma [[x^{1/2}, x^{-1/2}]]$ as follows: For $\alpha \in \mathfrak{h}$, $n \in \mathbb{N} + 1/2$, and $u \in V_{fer}$, let
\begin{multline}\label{define-sigma-twisted}
Y^\sigma (\alpha(-n) u,x) = Y^\sigma(\alpha_{-n -1/2} u,x) = \mathrm{Res}_{x_1} \mathrm{Res}_{x_0} \left( \frac{x_1 - x_0}{x} \right)^{1/2} x_0^{-n-1/2}\\
\cdot
\left(  x^{-1}_0\delta\left(\frac{x_1-x}{x_0}\right)
\alpha (x_1)^\sigma Y^\sigma(u,x)
-  (-1)^{|u|} x^{-1}_0\delta\left(\frac{x-x_1}{-x_0}\right) Y^\sigma
(u,x) \alpha (x_1)^\sigma  \right) .
\end{multline}
Then since $V_{fer} = \langle \alpha^{(j)}(-1/2) \mathbf{1} \; | \; j = 1,\dots, d \rangle$, equation (\ref{define-sigma-twisted}) defines $Y^\sigma(v,x)$ iteratively for any $v \in V_{fer}$. 

Recalling that the Virasoro element, $\omega_{fer}$, for the free fermionic VOSA $V_{fer}$ is given by (\ref{fermionic-Virasoro}), we have 
\begin{equation}
Y^\sigma(\omega_{fer}, x) = \frac{1}{2} \sum_{j = 1}^d Y^\sigma( \alpha^{(j)} (-1/2)_{-2} \alpha^{(j)}(-1/2)\mathbf{1},x)  = \sum_{n \in \mathbb{Z}}  L^\sigma(n) x^{-n-2}, 
\end{equation}
and thus
\begin{equation}\label{L-sigma-twisted}
L^\sigma(m) = \sum_{j = 1}^d \sum_{{n \in \mathbb{Z}}\atop{n > -\frac{n}{2}}} \left(n + \frac{m}{2}\right) \overline{\alpha^{(j)} (-n)} \, \overline{\alpha^{(j)} (n +m) }  + \frac{d}{16} \delta_{m,0}.
\end{equation}
Therefore
\begin{eqnarray*}
\lefteqn{ \left[L^\sigma(-1), Y^\sigma( \alpha^{(j)} (-1/2)\mathbf{1},x )\right] \ = \ \left[L^\sigma(-1), \alpha^{(j)} (x)^\sigma \right] }\\
&=& \sum_{r = 1}^d \sum_{m \in \mathbb{Z}_+} \sum_{n \in \mathbb{Z}} \left(m-\frac{1}{2} \right) \left[\overline{\alpha^{(r)} (-m)} \,  \overline{\alpha^{(r)} (m - 1)} , \, \overline{\alpha^{(j)} (n) }\right] x^{-n - 1/2} \\
&=& \sum_{r = 1}^d \sum_{m \in \mathbb{Z}_+} \sum_{n \in \mathbb{Z}} \left(m-\frac{1}{2} \right) \left( \overline{\alpha^{(r)} (-m)}  \left[\overline{\alpha^{(r)} (m - 1)}, \, \overline{\alpha^{(j)} (n)} \right] \right. \\
& & \quad \left.  -  \left[\overline{\alpha^{(j)} (n)},\, \overline{\alpha^{(r)} (-m )} \right] \overline{\alpha^{(r)} (m-1) }  \right) x^{-n - 1/2} \\
&=&  \sum_{m \in \mathbb{Z}_+}  \left( m-\frac{1}{2} \right) \left( \overline{\alpha^{(j)} (-m)} x^{m - 3/2}  - \overline{\alpha^{(j)} (m-1) } x^{-m - 1/2} \right)\\
&=& \sum_{m \in \mathbb{Z}}  \left( -m- \frac{1}{2} \right)  \overline{\alpha^{(j)} (m) }x^{-m - 3/2} \ = \ \frac{d}{d x} \alpha^{(j)}(x)^\sigma \\
&=& \frac{d}{dx} Y^\sigma (\alpha^{(j)}(-1/2)\mathbf{1}, x).
\end{eqnarray*}
It follows from \cite{Li-twisted}, that $M_\sigma$ is a weak $\sigma$-twisted module for $V_{fer}$.  It is also admissible.  

In \cite{FFR}, if $d$ is even, $M_\sigma$ is denoted by $CM(\mathbb{Z})$.  

By \cite{Li-twisted} as well as \cite{DZ}, in the case that $d = \mathrm{dim} \, \mathfrak{h}$ is even, $M_\sigma$ is irreducible and is the only irreducible admissible $\sigma$-twisted module for $V_{fer}$, up to isomorphism.  It is in fact also an ordinary $\sigma$-twisted $V_{fer}$-module, as we will see below when we discuss the $L^\sigma(0)$ grading.

In the case that $d$ is odd, $M_\sigma$ reduces as the direct sum of two irreducible admissible $\sigma$-twisted modules, and these two irreducibles are the only irreducible admissible $\sigma$-twisted  modules for $V_{fer}$, up to isomorphism.   In this case, setting
\[ W = \mbox{$\bigwedge$} \left[\overline{\beta^{(j)}_-(-m)} \mathbf{1}, \, \overline{\beta^{(j)}_+(-n)}\mathbf{1}, \,  \overline{\epsilon (-m) }\mathbf{1} \, \Big| \, \mbox{for $m \in \mathbb{Z}_+$, $n \in \mathbb{N}$, and $j = 1, \dots, l$} \right] \]
and letting $W = W^0 \oplus W^1$ be the decomposition of $W$ into even and odd subspaces, 
these two irreducibles are given by
\begin{eqnarray}
M_\sigma^\pm = \left(1 \pm  \, \overline{\epsilon(0)}\right) W^0 \oplus \left(1 \mp \, \overline{\epsilon(0)} \right) W^1,
\end{eqnarray}
and we have $M_\sigma = M_\sigma^- \oplus M_\sigma^+$.  The $M_\sigma^\pm$ are in fact ordinary $\sigma$-twisted modules for $V_{fer}$, as we shall see now by discussing the $L^\sigma(0)$ grading.

In terms of the polarization of $\mathfrak{h}$ with respect to the basis $\alpha^{(j)}$, we have from equation (\ref{L-sigma-twisted})
\begin{equation} 
L^{\sigma}(0)  =  \sum_{j = 1}^l  \Biggl( \sum_{m  \in \mathbb{Z}_+} \left( m  \overline{\beta^{(j)}_+( -m) } \, \overline{\beta_-^{(j)} (m) } + m  \overline{\beta^{(j)}_-( -m) } \, \overline{\beta_+^{(j)} (m) } \right)  \Biggr)  + L' + \frac{d}{16}
\end{equation}
where
\begin{equation}
L' = \left\{ \begin{array}{ll} 
0 & \mbox{if $d = 2l$}\\
\displaystyle{\frac{1}{2}\sum_{m\in \mathbb{Z}_+} m \overline{\epsilon (-m)} \, \overline{\epsilon(m) }  }& \mbox{if $d=2l+1$}
\end{array} \right. .
\end{equation}
Thus for $j = 1,\dots, l$, and $m \in \mathbb{Z}_+$, the $L^\sigma(0)$ grading is given by 
\begin{equation}
\mathrm{wt} \, \mathbf{1} =  \mathrm{wt} \, \overline{\beta^{(j)}_+ (0) }\mathbf{1} = \frac{d}{16},  \quad \mathrm{and} \quad  \mathrm{wt} \, \overline{\beta^{(j)}_\pm (-m)} \mathbf{1} = m + \frac{d}{16}, 
\end{equation}
for $d = 2l$, and if $d$ is odd, we also have 
\begin{equation}
 \mathrm{wt} \,  \overline{\epsilon(0) } \mathbf{1} =  \frac{d}{16}, \quad  \mathrm{and} \quad  \mathrm{wt} \, \overline{\epsilon (-m) } \mathbf{1} =  m +\frac{d}{16}.
\end{equation}

Therefore, for $d$ even, the graded dimension of $M_\sigma$ is
\begin{eqnarray}
\mathrm{dim}_q M_\sigma \! \! &=& \! \! q^{-c/24} \sum_{\lambda \in \mathbb{C}} (M_\sigma)_\lambda q^\lambda  \ = \ q^{-d/48} q^{d/16} 2^{d/2} \prod_{n\in \mathbb{Z}_+} (1 + q^n)^d \\
&=&  \! \! \mathfrak{f}_2 (q)^d \nonumber
\end{eqnarray} 
where $\mathfrak{f}_2$ is a classical Weber function as discussed in Remark \ref{Weber-remark}.  
For $d$ odd, the graded dimension of $M_\sigma$ is 
\begin{equation}
\mathrm{dim}_q M_\sigma \ =\  q^{-d/48} q^{d/16} 2^{(d+1)/2} \prod_{n\in \mathbb{Z}_+} (1 + q^n)^d \ = \ \sqrt{2} \mathfrak{f}_2 (q)^d ,
\end{equation} 
and the grading of each of the two submodules $M^\pm_\sigma$ is exactly half that of the graded dimension of $M_\sigma$.

\subsection{$\sigma$-twisted modules for free and lattice N=1  and N=2 VOSAs---the Ramond sectors}\label{N=1,2-Ramond-section}

Setting $M = V_{bos} \otimes M_\sigma$, we have that $M$ is a $\sigma$-twisted module for the N=1 VOSA, $V = V_{bos} \otimes V_{fer}$, and thus is naturally a representation of the Ramond algebra.  Specifically, we have that the N=1 Ramond algebra representation is given by (\ref{L-sigma-twisted}) and 
\begin{equation}
G^\sigma (n)  = \sum_{j = 1}^d  \sum_{m \in \mathbb{Z}} \alpha^{(j)}( m) \overline{\alpha^{(j)} (n- m)}
\end{equation}
with central charge $c = 3d/2$.

The graded trace for the $\sigma$-twisted module $M = V_{bos} \otimes M_\sigma$, i.e. the Ramond sector, for the N=1 VOSA, $V = V_{bos}\otimes V_{fer}$, for $d$ even, is
\begin{equation}
\quad \mathrm{dim}_q M   \ = \  \sqrt{2}^{d} \prod_{n\in \mathbb{Z}_+} (1 + q^n)^d (1-q^n)^{-d} \ = \    \left(\frac{\mathfrak{f}_2 (q)}{\eta(q)} \right)^d \  = \   \sqrt{2}^{d} \left(\frac{\eta (q^2)}{\eta(q)^2} \right)^d ,
\end{equation} 
and for $d$ odd is 
\begin{equation}
\quad \mathrm{dim}_q M  \ = \  \sqrt{2}  \left(\frac{\mathfrak{f}_2 (q)}{\eta(q)} \right)^d \  = \   \sqrt{2}^{d + 1} \left(\frac{\eta (q^2)}{\eta(q)^2} \right)^d .
\end{equation}

\begin{rema}{\em
Since 
\begin{equation}
\mathfrak{f} (q)\mathfrak{f}_1(q) \mathfrak{f}_2(q) = \sqrt{2}
\end{equation}
we have that 
\begin{equation}
\left( \mathrm{dim}_q V \right) \left(\mathrm{sdim}_q V \right) \left( \mathrm{dim}_q M \right) =  C_d \frac{\sqrt{2}^d}{\eta(q)^{3d}} ,
\end{equation}
where 
\begin{equation}\label{define-coefficient}
C_d = \left\{ \begin{array}{ll}  1 & \mbox{for $d$ even} \\
\sqrt{2}  & \mbox{for $d$ odd} 
\end{array} \right. .
\end{equation}

}
\end{rema}

For the free N=2 VOSA of central charge $c = 3d$ given by $V\otimes V$, where $V= V_{bos} \otimes V_{fer}$ is the $d$ free boson-fermion N=1 VOSA, we have that 
\begin{equation}
V_{bos} \otimes M_{\sigma} \otimes V_{bos} \otimes M_{\sigma}
\end{equation}
is a Ramond twisted sector, i.e. a $\sigma$-twisted $V \otimes V$-modules, 
where $M_\sigma$ is the Ramond twisted sector for the $d$ free fermion VOSA.  The graded dimension of this module is then of course the $q$-dimension of  $V_{bos} \otimes M_\sigma$ squared. 

The Ramond twisted module $(V_{bos} \otimes M_{\sigma})^{\otimes 2}$ for the N=2 VOSA $(V_{bos} \otimes V_{fer})^{\otimes 2}$ has $p,q$-dimension given by 
\begin{equation}
\mathrm{dim}_{p,q} V_{bos} \otimes M_{\sigma} \otimes V_{bos} \otimes M_{\sigma} = \eta(q)^{-2d} 2^{d} C_d^2 \prod_{n \in \mathbb{Z}_+} (1 + pq^n)^d(1 + p^{-1} q^n)^d .
\end{equation}

The N=2 Ramond algebra representation with central charge $3d$ is generated by
\begin{eqnarray}
\ \ \ \ \ \ G^{(1),\sigma} (n)  \! \! \! &=& \! \! \! \sum_{j = 1}^d  \sum_{m \in \mathbb{Z}} \left( \alpha^{(j)}( m) \overline{\alpha^{(j)} (n- m)} \otimes 1 + 1 \otimes \alpha^{(j)}( m) \overline{\alpha^{(j)} (n- m)} \right)\\
G^{(2), \sigma}(n) \! \! \! &=&  \! \! \! \sum_{j = 1}^d  \sum_{m \in \mathbb{Z}} \left( \alpha^{(j)}( m) \otimes \overline{\alpha^{(j)} (n- m)}  + \alpha^{(j)}( m) \otimes \overline{\alpha^{(j)} (n- m)} \right).
\end{eqnarray}

Replacing $V_{bos}$ with $V_L$ for $L$ a positive definite even lattice of rank $d$, and where $\mathfrak{h} = L \otimes_{\mathbb Z} \mathbb{C}$, we obtain the Ramond twisted sectors for the corresponding lattice N=1 and N=2 VOSAs.   That is $V_L \otimes M_\sigma$ is a $\sigma$-twisted $V_L \otimes V_{fer}$-module.  And of course the corresponding graded traces are multiplied by one or two factors of $\Theta(L)$, for N=1 and N=2 respectively, where $\Theta(L)$ is  the theta-function associated to the lattice.  If $L$ is a positive definite integral lattice, but not even, then $V_L \otimes M_\sigma$ is a $Id_V \otimes \sigma$-twisted $V_L \otimes V_{fer}$-module and again gives a representation of the N=1 Ramond algebra (or the N=2 Ramond algebra if tensored with itself).  But the parity-twisted modules for $V_L \otimes V_{fer}$, that is the $\sigma \otimes \sigma$-modules for each $\sigma$ the parity automorphism on the respective VOSA, have not yet been constructed.

\section{Mirror maps and mirror-twisted sectors for free N=2 VOSAs and extensions to lattice N=2 VOSAs}

In this section, we define two distinct mirror maps for free and lattice N=2 VOSAs.  These two mirror maps give distinct equivalence classes in the group of automorphisms for the VOSAs.  Thus these two mirror maps result in distinct mirror-twisted module structures.  We explicitly construct mirror-twisted modules for free N=2 VOSAs for one of the mirror maps, and also indicate how to construct mirror-twisted modules for lattice N=2 VOSAs for this mirror map.   We discuss the obstructions and state of knowledge for constructing mirror-twisted modules for the other mirror map.

\subsection{Two distinct mirror maps for free and lattice N=2 VOSAs}

For the free N=2 VOSA of central charge $c = 3d$ constructed in Section \ref{N=2-free-section}, and denoted $V \otimes V$ where $V$ is the N=1 VOSA of central charge $3d/2$, we can define a mirror map $\kappa$ as follows: let
\begin{equation}
\kappa: \alpha^\pm (-n)\mathbf{1} \mapsto \alpha^\mp (-n) \mathbf{1}
\end{equation}
for $n \in \frac{1}{2} \mathbb{Z}_+$, and $\alpha^\pm = \frac{1}{\sqrt{2}} (\alpha_{(1)} \mp i \alpha_{(2)})$ with $\alpha_{(1)} = \alpha \otimes \mathbf{1}$, and $\alpha_{(2)} = \mathbf{1} \otimes \alpha$, for $\alpha \in \mathfrak{h}$, and then extend to $V\otimes V$ via these generators.  This is equivalent to 
\begin{equation}\label{define-kappa} 
\kappa: \alpha_{(1)} (-n) \mathbf{1} \mapsto \alpha_{(1)} (-n) \mathbf{1} \quad \mbox{and} \quad \alpha_{(2)} (-n) \mathbf{1} \mapsto - \alpha_{(2)}(-n) \mathbf{1} ,
\end{equation}
extended to $V\otimes V$.
That is, $\kappa$ is the identity on the first tensor factor of $V\otimes V$, and acts as $-1$ on the generators of the second tensor factor.  Note then that $\kappa$ is the parity map on the second fermionic factor $\mathbf{1} \otimes V_{fer} = \langle \alpha^{(j)}_{(2)} (-1/2) \mathbf{1} \, | \, j = 1,\dots, d \rangle$.

Furthermore, if we let $L$ be a positive definite lattice, and $V_L$ the corresponding VOSA, then letting $\kappa$ be the lattice isometry $\alpha \mapsto - \alpha$, for $\alpha \in L$, we have that $\kappa$ lifts (not necessarily uniquely) to a VOSA automorphism on $V_L$.  Then this VOSA automorphism along with the parity map on $V_{fer}$ defines a mirror map, which we also denote by $\kappa$.  

Then note that $\kappa(\mu) = -\mu$, $\kappa (\tau^\pm) = \tau^\mp$ and $\kappa (\omega) = \omega$ for both the free N=2 VOSA and the lattice N=2 VOSA.  

Note however, that there is another mirror map on free and lattice N=2 VOSAs.  Letting $V_b$ denote either $V_{bos}$ or $V_L$, then we have the following mirror map on $(V_b \otimes V_{fer})^{\otimes 2}$: 
\begin{eqnarray}
\tilde{\kappa}:  (V_b \otimes V_{fer})\otimes (V_b \otimes V_{fer}) &\longrightarrow &  (V_b \otimes V_{fer})\otimes (V_b \otimes V_{fer})  \\
u \otimes v &\mapsto & (-1)^{|u||v| } v \otimes u,
\end{eqnarray}
for $u, v \in V_b \otimes V_{fer}$ of homogeneous sign.  
That is $\tilde{\kappa}$ is a signed permutation map for $(V_b \otimes V_{fer})^{\otimes 2}$.  And we have $\tilde{\kappa}(\mu) = -\mu$, $\tilde{\kappa} (\tau^\pm) = \tau^\mp$ and $\tilde{\kappa} (\omega) = \omega$ for both the free N=2 VOSA and the lattice N=2 VOSA.  

These two different mirror maps, $\kappa$ and $\tilde{\kappa}$, for free and lattice N=2 VOSAs, partition the free and lattice N=2 VOSAs into different eigenspaces.  In particular, 
they correspond to distinct conjugacy classes in the automorphism groups of the VOSAs and thus necessarily result in different mirror-twisted module structures.   Below we construct the $\kappa$-twisted modules for the free and lattice N=2 VOSAs.  The construction of the $\tilde{\kappa}$-twisted modules involves extending the construction of permutation twisted modules for the tensor product of a VOA with itself, as achieved by the author along with Dong and Mason in \cite{BDM}, to VOSAs.  This construction will be developed in a subsequent paper.  However, we make note of the following:

\begin{lem}  The mirror-twisted module structures for the free and lattice N=2 vertex operator superalgebras constructed using the $\kappa$ mirror map are not isomorphic to the mirror-twisted module structures constructed using the $\tilde{\kappa}$ mirror map. 
\end{lem}

\subsection{Mirror-twisted modules for free N=2 VOSAs for the mirror map $\kappa$}

To construct a $\kappa$-twisted module for the free N=2 VOSA, $V\otimes V$, where $V = V_{bos} \otimes V_{fer}$ ($d$ free bosons coupled with $d$ free fermions), we first construct a $\kappa$-twisted module for the $d$ free bosons, $V_{bos}$.  We will denote this $\kappa$-twisted $V_{bos}$-module by $M_\kappa$.  Then $V \otimes M_\kappa \otimes M_\sigma$ will be a $\kappa$-twisted module for $V \otimes V$.  

To construct the $\kappa$-twisted $V_{bos}$-module, $M_\kappa$, we first let $t$ again be a formal commuting variable, and form the affine Lie algebra
\[\hat{\mathfrak{h}}^b[\kappa] = \mathfrak{h} \otimes t^{1/2} \mathbb{C}[t,t^{-1}]\oplus \mathbb{C} \mathbf{k} ,\]
with Lie bracket relations
\begin{eqnarray}
\bigl[\mathbf{k}, \hat{\mathfrak{h}}^b[\kappa] \bigr] &=& 0 \label{commutation1B-twisted}\\
\bigl[\alpha \otimes t^m, \beta \otimes t^n \bigr] &=& m  \langle \alpha, \beta \rangle  \delta_{m + n,0} \mathbf{k} \label{commutation2B-twisted}
\end{eqnarray}
for $\alpha, \beta \in \mathfrak{h}$ and $m, n \in \mathbb{Z} +  \frac{1}{2}$, where we have extended the nondegenerate symmetric bilinear form on $\mathfrak{h}$.

Then $\hat{\mathfrak{h}}^b[\kappa]$ is a $((\mathbb{Z}+ \frac{1}{2}) \cup \{0\})$-graded Lie algebra 
\[\hat{ \mathfrak{h} }^b[\kappa] = \coprod_{n \in(\mathbb{Z}+  \frac{1}{2}) \cup \{0\} }\hat{ \mathfrak{h}}^b[\kappa]_n \]
where $\hat{\mathfrak{h}}^b[\kappa]_0 = \mathbb{C} \mathbf{k}$, and  $\hat{\mathfrak{h}}^b[\kappa]_n = \mathfrak{h} \otimes t^{-n}$ for $n \in \mathbb{Z}+  \frac{1}{2}$.  And $\hat{\mathfrak{h}}^b[\kappa]$ is a Heisenberg algebra with graded subalgebras 
\[\hat{\mathfrak{h}}^b [\kappa]_+ = \mathfrak{h} \otimes t^{-1/2} \mathbb{C}[t^{-1}] \qquad \mbox{and} \qquad \hat{\mathfrak{h}}^b [\kappa]_- = \mathfrak{h} \otimes t^{1/2} \mathbb{C}[t], \]
and we have $\hat{\mathfrak{h}}^b[\kappa] = \hat{\mathfrak{h}}^b[\kappa]_- \oplus  \hat{\mathfrak{h}}^b[\kappa]_+ \oplus \mathbb{C}\mathbf{k}$. 

Let $\mathbb{C}$ be the $(\hat{\mathfrak{h}}^b [\kappa]_- \oplus  \mathbb{C} \mathbf{k})$-module such that $\hat{\mathfrak{h}}^b[\kappa]_-$ acts trivially and $\mathbf{k}$ acts as $1$.  Set
\begin{equation}
M_\kappa =  U(\hat{\mathfrak{h}}^b[\kappa]) \otimes_{U(\hat{\mathfrak{h}}^b [\kappa]_-  \oplus \mathbb{C} \mathbf{k})} \mathbb{C} \cong S( \hat{\mathfrak{h}}^b[\kappa]_+),
\end{equation}
so that $M_\kappa$ is naturally isomorphic to the symmetric algebra of polynomials in $\hat{\mathfrak{h}}^b[\kappa]_+$; see Remark \ref{twisted-kappa-remark}.  It is also the universal enveloping algebra for $\hat{\mathfrak{h}}^b[\kappa]_+$.   Let $\alpha \in \mathfrak{h}$ and $n \in\mathbb{Z}+  \frac{1}{2}$.  We will use the notation
\[ \overline{\alpha(n) }= \alpha \otimes t^n \in \hat{\mathfrak{h}}^b[\kappa]\]
where the overline is meant to distinguish elements of $\hat{\mathfrak{h}}^b[\kappa]$ from elements of $\hat{\mathfrak{h}}^f$, used to construct the free fermionic theory.

Note that $M_\kappa$ is a $\hat{\mathfrak{h}}^b[\kappa]$-module with action induced from the commutation relations (\ref{commutation1B-twisted}) and (\ref{commutation2B-twisted}) given by 
\begin{eqnarray}
\mathbf{k} \overline{\beta (-m)} \mathbf{1}&=& \overline{\beta (-m)} \mathbf{1} \label{kappa-twisted-relations1}\\
\overline{\alpha  (n)} \, \overline{\beta (-m) } \mathbf{1} &=& \langle \alpha, \beta \rangle  n \delta_{m,n} \mathbf{1} \\
\overline{\alpha (-n)} \, \overline{\beta (-m)} \mathbf{1} &=& \overline{\beta(-m)} \, \overline{\alpha(-n)} \mathbf{1}
\end{eqnarray}
for either $\alpha, \beta \in \mathfrak{h}$ and $m, n \in \mathbb{N}+ \frac{1}{2}$, and where here $\mathbf{1} = \mathbf{1}_{M_\kappa} = 1$.

\begin{rema}\label{twisted-kappa-remark}
{\em Let $\{\alpha^{(1)}, \alpha^{(2)}, \dots, \alpha^{(d)}\}$ be an orthonormal basis for $\mathfrak{h}$.  Then setting $b^{(j)}_n = \overline{\alpha^{(j)}(-n)}\mathbf{1} $, for $n \in  \mathbb{N}+ \frac{1}{2}$,  we have 
\begin{eqnarray*}
M_\kappa &=& \mathbb{C} [\overline{\alpha^{(j)}(-n)} \mathbf{1}\, | \, \mbox{for $n \in \mathbb{N}+ \frac{1}{2}$ and $j = 1, \dots, d$} ] \\
&=& \mathbb{C} [b^{(1)}_{\frac{1}{2}}, b^{(1)}_{\frac{3}{2}}, \dots, b^{(2)}_{\frac{1}{2}}, b^{(2)}_{\frac{3}{2}}, \dots, b^{(d)}_{\frac{1}{2}}, b^{(d)}_{\frac{3}{2}},\dots ],
\end{eqnarray*}
and we have the following operators on $M_\kappa$ 
\begin{eqnarray}
\mathbf{k} &\mapsto&  1 \\
\overline{\alpha^{(j)} (n) } &\mapsto&  n  \frac{\partial}{\partial  b^{(j)}_n} \\
\overline{\alpha^{(j)}(-n) }& \mapsto& b^{(j)}_{n}, \label{kappa-multiplication}
\end{eqnarray}
for $j = 1, \dots, d$, $n \in  \mathbb{N}+ \frac{1}{2}$, and where the operators on the left of (\ref{kappa-multiplication}) are multiplication operators. 
}
\end{rema} 

Let $x$ be a formal commuting variable, and for $\alpha \in \mathfrak{h}$, set
\begin{equation} \alpha(x)^\kappa = \sum_{n \in \mathbb{Z}+ \frac{1}{2}} \overline{\alpha (n) } x^{-n-1} .
\end{equation}
Then for $j,k = 1, \dots, d$, we have 
\begin{eqnarray}
[\alpha^{(j)}(x_1)^\kappa, \alpha^{(k)}(x_2)^\kappa ] &=& \delta_{j,k} x_1^{1/2}x_2^{-1/2}\left( \frac{1}{(x_1 - x_2)^2} - \frac{1}{(-x_2 + x_1)^2} \right. \\
& & \quad \left. - \frac{1}{2} x_2^{-1} \left(  \frac{1}{(x_1 - x_2)} - \frac{1}{(-x_2 + x_1)} \right)  \right)\nonumber
\end{eqnarray}
implying that the $\alpha^{(j)}(x)^\kappa$, for $j = 1,\dots, d$, are mutually local.

For $v \in V_{bos}$, define $Y^\kappa(v,x) : M_\kappa \longrightarrow M_\kappa [[x^{1/2}, x^{-1/2}]]$ as follows: For $\alpha \in \mathfrak{h}$, $n \in \mathbb{Z}_+$, and $u \in V_{bos}$, let
\begin{multline}\label{define-kappa-twisted}
Y^\kappa (\alpha(-n) u,x) = Y^\kappa(\alpha_{-n} u,x) = \mathrm{Res}_{x_1} \mathrm{Res}_{x_0} \left( \frac{x_1 - x_0}{x} \right)^{1/2} x_0^{-n}\\
\cdot
\left(  x^{-1}_0\delta\left(\frac{x_1-x}{x_0}\right)
\alpha (x_1)^\kappa Y^\kappa(u,x)
-  x^{-1}_0\delta\left(\frac{x-x_1}{-x_0}\right) Y^\kappa
(u,x) \alpha (x_1)^\kappa  \right) .
\end{multline}
Then since $V_{bos} = \langle \alpha^{(j)}(-1) \mathbf{1} \; | \; j = 1,\dots, d \rangle$, equation (\ref{define-kappa-twisted}) defines $Y^\kappa (v,x)$ recursively for any $v \in V_{bos}$. 

From (\ref{define-kappa-twisted}) and $Y^\kappa(\omega_{bos}, x) = \frac{1}{2} \sum_{j=1}^d Y^\kappa( \alpha^{(j)}_{-1} \alpha^{(j)} (-1) \mathbf{1}, x) = \sum_{n \in \mathbb{Z}} L^\kappa(n) x^{-n-2}$, we have that 
\begin{equation}\label{L-kappa-twisted}
L^\kappa(m) = \sum_{j=1}^d \sum_{{n \in \mathbb{Z}}\atop{n \geq - \frac{n}{2}}}  \overline{\alpha^{(j)} (-n - 1/2)} \, \overline{\alpha^{(j)} (m+n +1/2) }  + \tilde{L}^\kappa (n) + \frac{d}{16} \delta_{m,0},
\end{equation}
where
\begin{equation}
\tilde{L}^\kappa (n) = \left\{ \begin{array}{ll}
0 & \mbox{if $n$ is even}\\
\frac{1}{2} \sum_{j=1}^d \overline{\alpha^{(j)} (n/2)} \, \overline{\alpha^{(j)} (n/2) } & \mbox{if $n$ is odd}
\end{array}
\right. .
\end{equation}
Thus
\begin{eqnarray*}
\lefteqn{ \left[L^\kappa(-1), Y^\kappa( \alpha^{(j)} (-1)\mathbf{1},x) \right] \ = \ \left[L^\kappa(-1), \alpha^{(j)} (x)^\kappa \right] }\\
&=& \sum_{r = 1}^d  \sum_{n \in  \mathbb{Z}+ \frac{1}{2}}  \biggl(\sum_{m \in  \mathbb{Z}_++ \frac{1}{2}}\left[\overline{\alpha^{(r)} (-m)} \,  \overline{\alpha^{(r)} (m - 1)} , \, \overline{\alpha^{(j)} (n) }\right]  \\
& & \quad  + \frac{1}{2} \left[ \overline{\alpha^{(r)} (-1/2)}^2, \overline{\alpha^{(j)} (n)} \right] \biggr) x^{-n - 1} \\
&=& \sum_{r = 1}^d \sum_{n \in 
\mathbb{Z}+ \frac{1}{2}}  \biggl( \sum_{m \in  \mathbb{Z}_++ \frac{1}{2}} \biggl( \overline{\alpha^{(r)} (-m)}  \left[\overline{\alpha^{(r)} (m - 1)}, \, 
\overline{\alpha^{(j)} (n)} \right]  \\
& & \quad  -  \left[\overline{\alpha^{(j)} (n)},\, \overline{\alpha^{(r)} (-m )} \right] \overline{\alpha^{(r)} (m-1) }  \biggr) + \frac{1}{2} \overline{\alpha^{(r)} (-1/2)} \left[ \overline{\alpha^{(r)} (-1/2)}, \overline{\alpha^{(j)} (n)} \right] \\
& & \quad - \frac{1}{2} \left[ \overline{\alpha^{(j)} (n)}, \overline{\alpha^{(r)} (-1/2)} \right] \overline{\alpha^{(r)} (-1/2)}  \biggr) x^{-n - 1} \\
&=&  \sum_{m \in  \mathbb{Z}_+ + \frac{1}{2}} \left( (m-1) \overline{\alpha^{(j)} (-m)} x^{m - 2}  - m \overline{\alpha^{(j)} (m-1) } x^{-m - 1} \right) - \frac{1}{4} \overline{\alpha^{(j)} (-1/2)} x^{-3/2} \\
& & \quad - \frac{1}{4} \overline{\alpha^{(j)} (-1/2)} x^{-3/2}\\
&=& \sum_{m \in  \mathbb{Z} + \frac{1}{2}}  (-m-1)  \overline{\alpha^{(j)} (m) }x^{-m - 2} \ = \ \frac{d}{d x} \alpha^{(j)}(x)^\kappa \\
&=& \frac{d}{d x} Y^\kappa (\alpha^{(j)}(-1)\mathbf{1}, x).
\end{eqnarray*}
It follows from \cite{Li-twisted}, that $M_\kappa$ is a weak $\kappa$-twisted module for $V_{bos}$.  

Note that from (\ref{L-kappa-twisted}), we have that the $L^\kappa(0)$-grading of $M_\kappa$ implies that $M_\kappa$ is an ordinary $\kappa$-twisted $V_{bos}$-module with graded $q$-dimension given by 
\begin{equation}
\mathrm{dim}_q M_\kappa \ = \ q^{-d/24} q^{d/16} \prod_{n \in \mathbb{Z}_+} (1 + q^{n/2}) \ = \ \sqrt{2}^{-d} \mathfrak{f}_2 ( q^{1/2})^d \ = \ \left( \frac{\eta(q)}{\eta( q^{1/2})} \right)^d. 
\end{equation}

Setting $Y^\kappa( u \otimes v \otimes w, x) = Y(u,x)\otimes Y^\kappa( v,x) \otimes Y^\sigma(w, x)$, for $u \in V$, $v \in V_{bos}$ and $w \in V_{fer}$, we have that $V \otimes M_\kappa \otimes M_\sigma$ is a $\kappa$-twisted module for $V \otimes V = V \otimes V_{bos} \otimes V_{fer}$, where
\begin{eqnarray}
Y^\kappa(v_{(1)} , x) &=& Y^\kappa(v \otimes \mathbf{1}, x) \ = \ Y(v, x) \otimes Id_V\\
Y^\kappa( \alpha_{(2)}(-1)\mathbf{1}, x) &=& \sum_{n \in \mathbb{Z} + \frac{1}{2}} \overline{\alpha_{(2)} (n)}   x^{-n-1}   \\
Y^\kappa ( \alpha_{(2)}(-1/2)\mathbf{1}, x) &=& \sum_{n \in \mathbb{Z}} \overline{\alpha_{(2)}(n)} x^{-n-1}.
\end{eqnarray}

Then the $\kappa$-twisted $V\otimes V$-module, $V \otimes M_\kappa \otimes M_\sigma$ is in fact an ordinary $\kappa$-twisted module with $q$-dimension
\begin{eqnarray}
\ \ \ \ \ \ \mathrm{dim}_q (V \otimes M_\kappa \otimes M_\sigma) \! \! \!  &=&  \! \! \!  C_d \left(\frac{\sqrt{2}^{-1}\mathfrak{f}(q) \mathfrak{f}_2 (q)  \mathfrak{f}_2(q^{1/2})}{\eta(q)} \right)^d  =  C_d \left( \frac{\mathfrak{f}_2(q^{1/2})}{\eta(q^{1/2})} \right)^d \\
&=& \! \!  \! C_d \sqrt{2}^{d}  \left(\frac{\eta(q) }{\eta(q^{1/2})^2} \right)^d ,\nonumber
\end{eqnarray} 
where $C_d$ is given by (\ref{define-coefficient}).
Of course since $\kappa (\mu) = -\mu$, there is no zero mode for the mirror-twisted vertex operator associated to $\mu$ and thus no notion of $p,q$-dimension. 

\begin{rema} {\em
Note the similarity between the $q$-dimension of the $\kappa$-twisted $V \otimes V$-module, $V \otimes M_\kappa \otimes M_\sigma$ and the $q$-dimension of the N=1 Ramond twisted sector $V_{bos} \otimes M_\sigma$ for $V$.  That is, we have that 
\begin{equation}
\mathrm{dim}_q (V_{bos} \otimes M_\sigma) = \mathrm{dim}_{q^2}  (V \otimes M_\kappa \otimes M_\sigma) .
\end{equation} 
It is not clear the reason or significance of this similarity, but such similarities have been noted before, as in for instance \cite{IK2009}.
}
\end{rema}

\begin{rema}{\em
Note that $V \otimes M_\kappa \otimes M_\sigma = V_{bos} \otimes V_{fer} \otimes M_\kappa \otimes M_\sigma$ naturally contains the subspace $V_{fer} \otimes M_\sigma$ which in the notation of \cite{FFR} is $CM(\mathbb{Z}  + \frac{1}{2}) \otimes CM( \mathbb{Z})$, and is a $(Id_{V_{fer}} \otimes \sigma)$-twisted $V_{fer} \otimes V_{fer}$-module. 
}
\end{rema}

\begin{rema}{\em
Let $L$ be a positive definite lattice and $\kappa: \alpha \mapsto - \alpha$, for $\alpha \in L$, a lattice isometry.  Following \cite{Lepowsky1985}, \cite{DL}, and \cite{Xu}, we can lift $\kappa$ to an order two automorphism of $V_L$ and form the $\kappa$-twisted $V_L$-module, denoted $V_L^T$.  Then $V_L \otimes V_{fer} \otimes V_L^T \otimes M_\sigma$ is a $\kappa$-twisted $(V_L \otimes V_{fer})^{\otimes 2}$-module.  Note that $\kappa$ restricted to the Heisenberg part of $V_L$ is $\kappa$ acting as minus one on $\mathfrak{h} = L \otimes_{\mathbb{Z}} \mathbb{C}$ as in (\ref{define-kappa}) for the second tensor factor. }
\end{rema}

\section{$\sigma_\xi$-twisted modules for free and lattice N=2 VOSAs}\label{sigma-twisted-section}

For the free N=2 VOSA, $V \otimes V$, of central charge $c = 3d$ constructed in Section \ref{N=2-free-section}, we have a $J(0)$-grading with $J(0)( \alpha^\pm (-n) \mathbf{1}) = 0$,  for $n \in \mathbb{Z}_+$, and $J(0) (\alpha^{\pm}(-r)\mathbf{1}) = \pm \alpha^\pm (-r) \mathbf{1}$ for $r \in \mathbb{N} + \frac{1}{2}$.   Thus we can extend the N=2 Neveu-Schwarz algebra automorphism $\sigma_\xi$ to a VOSA automorphism of  $V \otimes V$ as follows:  
\begin{equation}
\sigma_\xi : \quad  \alpha^\pm(-n) \mathbf{1}  \mapsto  \alpha^\pm(-n) \mathbf{1} \qquad   \alpha^\pm (-r) \mathbf{1} \mapsto \xi^{\pm 1} \alpha^\pm (-r) \mathbf{1} 
\end{equation}
or more generally $\sigma_\xi(v) = \xi^n v$ if $J(0) v = nv$, for $n \in \mathbb{Z}$.

If $\xi$ is a $k$-th root of unity for $k \in \mathbb{Z}_+$, then this VOSA automorphism $\sigma_\xi$ is of finite order, and we can consider the $\sigma_\xi$-twisted $V \otimes V$-modules.  Fix $\eta = e^{2 \pi i/k}$ for $k \geq 3$, and fix $\xi = \eta^j$ to be a primitive $k$-th root of unity for $1\leq j < k$.    (The case for $k = 2$ was already constructed in Section \ref{N=1,2-Ramond-section}.)   We will construct the $\sigma_\xi =\sigma_{\eta^j}$-twisted sectors by first constructing a $\sigma_\xi$-twisted module for $V_{fer} \otimes V_{fer}$.  

Consider the vector space $\mathfrak{h} \oplus \mathfrak{h}$ with the nondegenerate symmetric bilinear form on $\mathfrak{h}$ extended to $\mathfrak{h} \oplus \mathfrak{h}$ by $\langle (\alpha_1, \beta_1), (\alpha_2, \beta_2) \rangle = \langle \alpha_1, \alpha_2 \rangle + \langle \beta_1, \beta_2 \rangle$.   Define the following subspaces of $\mathfrak{h} \oplus \mathfrak{h}$,
\begin{equation}
\mathfrak{h}^\pm = \mathrm{span}_{\mathbb{C}} \left\{ \alpha^\pm = \frac{1}{\sqrt{2}} ((\alpha,0) \mp i (0,\alpha) )\, \Big| \, \mbox{for $\alpha \in \mathfrak{h}$}\right\} .
\end{equation}
Note that $\langle \alpha^+ , \beta^+ \rangle = \langle \alpha^-, \beta^- \rangle = 0$ for $\alpha^\pm, \beta^\pm \in \mathfrak{h}^\pm$, and $\langle \alpha^+, \beta^- \rangle = \langle \alpha, \beta \rangle$.

Form the affine Lie superalgebra
\begin{equation}
\widehat{(\mathfrak{h} \oplus \mathfrak{h})} {}^f[\sigma_\xi] = \left( \left(\mathfrak{h}^- \otimes t^{1/2 -j/k}\mathbb{C}[t,t^{-1}] \right) \oplus \left( \mathfrak{h}^+ \otimes t^{1/2 + j/k}\mathbb{C}[t,t^{-1}] \right) \right) \oplus \mathbb{C} \mathbf{k}
\end{equation}
with $\mathbb{Z}_2$-grading given by $\mathrm{sgn}(\alpha^\pm\otimes t^n)= 1$ for $n \in  \mathbb{Z} + \frac{1}{2} \pm \frac{j}{k}$, and $\mathrm{sgn}(\mathbf{k})= 0$, and Lie super-bracket relations given by 
\begin{eqnarray}\label{sigma_xi-supercommutation1}
\bigl[\mathbf{k}, \widehat{(\mathfrak{h} \oplus \mathfrak{h})} {}^f[\sigma_\xi] \bigr] = \bigl[\alpha^\pm \otimes t^m, \beta^\pm \otimes t^n \bigr] &=& 0 \\
\bigl[\alpha^+ \otimes t^m, \beta^- \otimes t^n \bigr] = \langle \alpha^+, \beta^- \rangle  \delta_{m + n,0} \mathbf{k} &=&  \langle \alpha, \beta \rangle \delta_{m+n,0} \mathbf{k}, \label{sigma_xi-supercommutation2}
\end{eqnarray}
for $\alpha^\pm, \beta^\pm \in \mathfrak{h}^\pm$ and $m, n \in \mathbb{Z} + \frac{1}{2} \pm \frac{j}{k}$.

Then $\widehat{(\mathfrak{h} \oplus \mathfrak{h})} {}^f[\sigma_\xi]$ is a $\left(\mathbb{Z} + \frac{1}{2} - \frac{j}{k}\right) \cup \left(\mathbb{Z} + \frac{1}{2} + \frac{j}{k}\right)$-graded Lie superalgebra 
\[\widehat{( \mathfrak{h} \oplus \mathfrak{h}) } {}^f[\sigma_\xi] = \coprod_{n \in \left(\mathbb{Z} + \frac{1}{2} - \frac{j}{k}\right) \cup \left(\mathbb{Z} + \frac{1}{2} + \frac{j}{k}\right)}\widehat{( \mathfrak{h} \oplus \mathfrak{h})} {}^f[\sigma_\xi]_n \]
where $\widehat{(\mathfrak{h} \oplus \mathfrak{h})} {}^f[\sigma_\xi]_0 =  \mathbb{C} \mathbf{k}$,  $\widehat{(\mathfrak{h} \oplus \mathfrak{h})} {}^f[\sigma_\xi]_n = \mathfrak{h}^+ \otimes t^{-n}$, for $n \in \mathbb{Z} + \frac{1}{2} - \frac{j}{k}$, and 
$\widehat{(\mathfrak{h} \oplus \mathfrak{h})} {}^f[\sigma_\xi]_n = \mathfrak{h}^- \otimes t^{-n}$, for $n \in \mathbb{Z} + \frac{1}{2} + \frac{j}{k}$.   Note that  $\widehat{(\mathfrak{h} \oplus \mathfrak{h})} {}^f[\sigma_\xi]$ is a Heisenberg superalgebra.  

Then $\widehat{(\mathfrak{h} \oplus \mathfrak{h})} {}^f [\sigma_\xi]$ has the following graded subalgebras
\begin{equation*}
\widehat{(\mathfrak{h} \oplus \mathfrak{h})} {}^f[\sigma_\xi]_+ = \! \!  \! \coprod_{ n \in \mathbb{Z} + \frac{1}{2} \pm \frac{j}{k}\atop{n<0}} \! \! \! \widehat{(\mathfrak{h} \oplus \mathfrak{h})}{}^f[\sigma_\xi]_n, \quad \mbox{and} \quad 
\widehat{(\mathfrak{h} \oplus \mathfrak{h})} {}^f[\sigma_\xi]_- = \! \! \! \coprod_{ n \in \mathbb{Z} + \frac{1}{2} \pm \frac{j}{k}\atop{n>0}}  \! \! \! \widehat{(\mathfrak{h} \oplus \mathfrak{h})}{}^f[\sigma_\xi]_n,
\end{equation*}
and we have $\widehat{(\mathfrak{h} \oplus \mathfrak{h})} {}^f [\sigma_\xi] = \widehat{(\mathfrak{h} \oplus \mathfrak{h})} {}^f [\sigma_\xi]_- \oplus \widehat{(\mathfrak{h} \oplus \mathfrak{h})} {}^f [\sigma_\xi]_+ \oplus \mathbb{C} \mathbf{k}$.

Let $\mathbb{C}$ be the $(\widehat{(\mathfrak{h} \oplus \mathfrak{h})} {}^f[\sigma_\xi]_- \oplus \mathbb{C} \mathbf{k})$-module such that  $\widehat{(\mathfrak{h} \oplus \mathfrak{h})} {}^f[\sigma_\xi]_-$ acts trivially and  $\mathbf{k}$ acts as 1.  Set
\begin{equation}
M_{\sigma_\xi}  = U\left(\widehat{(\mathfrak{h} \oplus \mathfrak{h})} {}^f[\sigma_\xi]\right) \otimes_{U\left(\widehat{(\mathfrak{h} \oplus \mathfrak{h})} {}^f[\sigma_\xi]_- \oplus \mathbb{C} \mathbf{k}\right)} \mathbb{C} \cong \mbox{$\bigwedge$}\left( \widehat{(\mathfrak{h} \oplus \mathfrak{h})} {}^f [\sigma_\xi]_+ \right) ,
\end{equation}
so that $M_{\sigma_\xi}$ is naturally isomorphic to the algebra of polynomials in the anticommuting elements of $ \widehat{(\mathfrak{h} \oplus \mathfrak{h})} {}^f [\sigma_\xi]_+$; see Remark \ref{twisted-NS-algebra-remark}. 

Let $\alpha^\pm \in \mathfrak{h}^\pm$ and $n \in \mathbb{Z} + \frac{1}{2} \pm \frac{j}{k}$, respectively.  We will use the notation
\[ \alpha^\pm (n) = \alpha^\pm \otimes t^n.\]

Then $M_{\sigma_\xi}$ is a $\widehat{(\mathfrak{h} \oplus \mathfrak{h})} {}^f[\sigma_\xi]$-module with action induced from the supercommutation relations (\ref{sigma_xi-supercommutation1}) and (\ref{sigma_xi-supercommutation2}) given by 
\begin{eqnarray}
\mathbf{k} \beta^\pm (-m) \mathbf{1} &=& \beta^\pm (-m) \mathbf{1} \\
\alpha^\pm (n) \beta^\pm (-m) \mathbf{1} &=& 0\\
\alpha^\mp (n') \beta^\pm (-m) \mathbf{1} &=& \langle \alpha, \beta \rangle \delta_{m,n'} \mathbf{1}\\
\alpha^\pm(-n'') \beta^\pm (-m) \mathbf{1} &=& - \beta^\pm (-m) \alpha^\pm (-n'') \mathbf{1}\\
\alpha^\mp (-n''') \beta^\pm(-m) \mathbf{1} &=&  - \beta^\pm(-m) \alpha^\mp (-n''') \mathbf{1}
\end{eqnarray}
for $\alpha^\pm, \beta^\pm \in \mathfrak{h}^\pm$, $m \in \mathbb{Z} + \frac{1}{2} \mp \frac{j}{k}$, $m>0$, respectively, 
$n, n''' \in \mathbb{Z} + \frac{1}{2} \pm \frac{j}{k}$, $n, n''' >0$, respectively, and $n',n''  \in \mathbb{Z} + \frac{1}{2} \mp \frac{j}{k}$, $n', n'' >0$, respectively.

\begin{rema}\label{twisted-NS-algebra-remark}{\em Again let $\{\alpha^{(1)}, \dots, \alpha^{(d)} \}$ be an orthonormal basis for $\mathfrak{h}$.  Let $a^{(m), \pm}_n$ for $n \in \mathbb{Z} + \frac{1}{2} \pm \frac{j}{k}$, $n >0$, respectively, be formal variables.  Then, for instance if $0 < \frac{j}{k} < \frac{1}{2}$, we have that $\widehat{(\mathfrak{h} \oplus \mathfrak{h})} {}^f [\sigma_\xi]$ acts on the space
\begin{multline}\label{NS-twisted-space}
\mbox{$\bigwedge$} \left[ a^{(1), +}_{\frac{1}{2} - \frac{j}{k}},  a^{(1), +}_{\frac{3}{2} - \frac{j}{k}},  a^{(1), +}_{\frac{5}{2} - \frac{j}{k}}, \dots, a^{(2), +}_{\frac{1}{2} - \frac{j}{k}},  a^{(2), +}_{\frac{3}{2} - \frac{j}{k}},  \dots, a^{(d), +}_{\frac{1}{2} - \frac{j}{k}},  a^{(d), +}_{\frac{3}{2} - \frac{j}{k}},  \dots,  \right. \\
\left. a^{(1), -}_{\frac{1}{2} + \frac{j}{k}},  a^{(1), -}_{\frac{3}{2} + \frac{j}{k}},  a^{(1), -}_{\frac{5}{2} + \frac{j}{k}}, \dots, a^{(2), -}_{\frac{1}{2} + \frac{j}{k}},  a^{(2), -}_{\frac{3}{2} + \frac{j}{k}},  \dots, a^{(d), -}_{\frac{1}{2} + \frac{j}{k}},  a^{(d), -}_{\frac{3}{2} + \frac{j}{k}},\dots\right]
\end{multline} 
by
\begin{eqnarray}
\mathbf{k} & \mapsto & 1 \\
\alpha^{(m), \pm} (n) & \mapsto & \frac{\partial \ \  \ }{\partial a^{(m), \mp}_n}\\
\alpha^{(m), \mp} (-n) & \mapsto & a^{(m), \mp}_n,
\end{eqnarray}
for $m = 1, \dots, d$, and $n \in \mathbb{N} + \frac{1}{2} \pm \frac{j}{k}$, respectively.  That is the space (\ref{NS-twisted-space}) is isomorphic to $M_{\sigma_\xi}$ as an $\widehat{(\mathfrak{h} \oplus \mathfrak{h})} {}^f[\sigma_\xi]$-module.  And similarly if $\frac{1}{2} < \frac{j}{k} < 1$.
}
\end{rema}

For $\alpha^\pm \in \mathfrak{h}$, set 
\begin{equation}
\alpha^\pm (x)^{\sigma_\xi}= \sum_{ n \in \mathbb{Z} + \frac{1}{2} \pm \frac{j}{k}} \alpha^\pm(n) x^{-n-\frac{1}{2}}.
\end{equation}
Then
\begin{eqnarray}
\left[\alpha^{(m), \pm} (x_1)^{\sigma_\xi}, \alpha^{(n), \pm}(x_2)^{\sigma_\xi}\right] \! \! \! \! &=& \! \! \! \! 0\\
\qquad \ \ \left[\alpha^{(m), +} (x_1)^{\sigma_\xi} , \alpha^{(n), -}(x_2)^{\sigma_\xi} \right] \! \! \! \!  &=& \! \! \!  \! \delta_{m,n} \, x_1^{1- \frac{j}{k}} x_2^{-1 + \frac{j}{k}} \! \left( \frac{1}{(x_1 - x_2)} - \frac{1}{(-x_2 + x_1)} \right)
\end{eqnarray}
for $m, n = 1, \dots, d$, implying that the $\alpha^{(m), \pm}(x)^{\sigma_\xi}$, for $m = 1,\dots, d$, are mutually local.   

For $v \in V_{fer} \otimes V_{fer}$, define $Y^{\sigma_\xi} (v,x) : M_{\sigma_\xi} \longrightarrow M_{\sigma_\xi} [[x^{1/k}, x^{-1/k}]]$ as follows: For $\alpha^\pm \in \mathfrak{h}^\pm$, $n \in \mathbb{N} + \frac{1}{2}$, and $u \in V_{fer} \otimes V_{fer}$, let
\begin{multline}\label{define-sigma-xi-twisted}
Y^{\sigma_\xi} (\alpha^\pm(-n) u,x) = Y^{\sigma_\xi} (\alpha^\pm_{-n - 1/2} u,x) = \mathrm{Res}_{x_1} \mathrm{Res}_{x_0} \left( \frac{x_1 - x_0}{x} \right)^{\pm 1/k} x_0^{-n- 1/2}\\
\cdot
\left(  x^{-1}_0\delta\left(\frac{x_1-x}{x_0}\right)
\alpha^\pm (x_1)^{\sigma_\xi} Y^{\sigma_\xi} (u,x)
-  (-1)^{|u|} x^{-1}_0\delta\left(\frac{x-x_1}{-x_0}\right) Y^{\sigma_\xi}
(u,x) \alpha^\pm (x_1)^{\sigma_\xi}  \right) .
\end{multline}
Then since $V_{fer} \otimes V_{fer}  = \langle \alpha^{(m), \pm} (-1/2) \mathbf{1} \; | \; m = 1,\dots, d \rangle$, equation (\ref{define-sigma-xi-twisted}) defines $Y^{\sigma_\xi} (v,x)$ recursively for all $v \in V_{fer}\otimes V_{fer}$. 

Since
\begin{multline*} \omega_{fer} \otimes \mathbf{1} + \mathbf{1} \otimes \omega_{fer} \\
= \frac{1}{2} \sum_{m=1}^d \left( \alpha^{(m), +} (-3/2) \alpha^{(m), -} (-1/2) \mathbf{1} + \alpha^{(m), - } (-3/2) \alpha^{(m), +} (-1/2)\mathbf{1}\right) ,
\end{multline*}
and 
\begin{multline}
\sum_{n \in \mathbb{Z}} L^{\sigma_\xi} (n) x^{-n-2} \\
= \frac{1}{2} \sum_{m = 1}^d \left( Y^{\sigma_\xi} ( \alpha^{(m), +}_{-2} \alpha^{(m),-} (-1/2) \mathbf{1}, x) + Y^{\sigma_\xi} ( \alpha^{(m), -}_{-2} \alpha^{(m),+}(-1/2) \mathbf{1}, x \right),
\end{multline} 
from (\ref{define-sigma-xi-twisted}) we have that 
\begin{multline}\label{L-sigma-xi}
L^{\sigma_\xi} (n) = \sum_{m=1}^d \Biggl(  \sum_{{r \in \mathbb{Z}+ \frac{1}{2} - \frac{j}{k}}\atop{r> - \frac{n}{2}}}  \Bigl(r + \frac{n}{2}\Bigr)  \alpha^{(m),+} (-r) \alpha^{(m), -} (r + n)  \\
+ 
 \sum_{{r \in \mathbb{Z}+ \frac{1}{2} + \frac{j}{k}}\atop{r> - \frac{n}{2} }}  \Bigl(r + \frac{n}{2}\Bigr)  \alpha^{(m),-} (-r) \alpha^{(m), +} (r + n)
\Biggr) + \frac{j^2 d}{2k^2} \delta_{n,0}.
\end{multline}

Thus
\begin{eqnarray*}
\lefteqn{\left[ L^{\sigma_\xi}(-1) , Y^{\sigma_\xi} (\alpha^{(m), \pm} (-1/2)\mathbf{1}, x)\right] \ =  \ \left[L^{\sigma_\xi}(-1), \alpha^{(m), \pm} (x)^{\sigma_\xi} \right] }\\
&=&   \sum_{n=1}^d  \sum_{ s \in \mathbb{Z} + \frac{1}{2} \pm \frac{j}{k}} \Biggl( \sum_{{r \in \mathbb{Z}+ \frac{1}{2} - \frac{j}{k}}\atop{r>\frac{1}{2}}}  \Bigl(r - \frac{1}{2}\Bigr)  \left[ \alpha^{(n),+} (-r ) \alpha^{(n), -} (r-1) , \alpha^{(m),\pm}(s) \right]  \\
& & + \sum_{{r \in \mathbb{N}+ \frac{1}{2} + \frac{j}{k}}\atop{r>\frac{1}{2}}} \Bigl(r - \frac{1}{2}\Bigr) \left[ \alpha^{(n),-} (-r ) \alpha^{(n), +} (r-1)  ,   \alpha^{(m),\pm}(s) \right] \Biggr) x^{-s-\frac{1}{2}} \\
&=&   \sum_{n=1}^d  \sum_{ s \in \mathbb{Z} + \frac{1}{2} \pm \frac{j}{k}} \Biggl( \sum_{{r \in \mathbb{Z}+ \frac{1}{2} - \frac{j}{k}}\atop{r>\frac{1}{2}}}  \Bigl(r - \frac{1}{2}\Bigr) \left(   \alpha^{(n),+} (-r) \left[\alpha^{(n), -} (r-1) , \alpha^{(m),\pm}(s) \right] \right. \\
& & \left. - \left[ \alpha^{(n),+} (-r ) , \alpha^{(m),\pm}(s) \right] \alpha^{(n), -} (r-1)  \right)  \\
& & + \sum_{{r \in \mathbb{Z}+ \frac{1}{2} + \frac{j}{k}}\atop{r>\frac{1}{2}}}  \Bigl(r - \frac{1}{2}\Bigr)  \left( \alpha^{(n),-} (-r ) \left[\alpha^{(n), +} (r-1)  ,   \alpha^{(m),\pm}(s) \right]  \right. \\
& &  \left. - \left[ \alpha^{(n),-} (-r) ,   \alpha^{(m),\pm}(s) \right]  \alpha^{(n), +} (r-1) \right) \Biggr) x^{-s-\frac{1}{2}} \\
&=&  \sum_{{r \in \mathbb{Z}+ \frac{1}{2} \mp \frac{j}{k}}\atop{r>\frac{1}{2}}} \Bigl(r - \frac{1}{2}\Bigr)  \alpha^{(m),\pm} (-r) x^{r-\frac{3}{2}} +  \sum_{{r \in \mathbb{Z} + \frac{1}{2} \pm \frac{j}{k}}\atop{r>\frac{1}{2}}} \Bigl(-r + \frac{1}{2}\Bigr) \alpha^{(m), \pm} (r-1) x^{-r - \frac{1}{2}}  \\
&=& \sum_{ r \in \mathbb{Z} + \frac{1}{2} \pm \frac{j}{k}}  \Bigl(-r -\frac{1}{2} \Bigr)   \alpha^{(m),\pm} (r) x^{-r-\frac{3}{2}} \ = \ \frac{d}{d x} \alpha^{(m), \pm} (x)^{\sigma_\xi} \ \ \ \ \ \ \ \ \ \ \ \\
&=& \frac{d}{d x} Y^{\sigma_\xi} (\alpha^{(m), \pm}(-1/2)\mathbf{1}, x).
\end{eqnarray*}
It follows from \cite{Li-twisted} that $M_{\sigma_\xi}$ is a $\sigma_\xi$-twisted module for $V_{fer} \otimes V_{fer}$.  

\begin{rema}{\em  From \cite{Li-twisted}, we have that $M_{\sigma_\xi}$ is the unique, up to isomorphism, irreducible $\sigma_\xi$-twisted $V_{fer}\otimes V_{fer}$-module.  }
\end{rema}

Since 
\begin{multline}
L^{\sigma_\xi} (0) =  \sum_{m = 1}^d \Biggl( \sum_{{r \in \mathbb{Z} + \frac{1}{2} - \frac{j}{k}}\atop{r>0}} r \  \alpha^{(m), +}(-r) \alpha^{(m), -} (r)\\ 
+ \sum_{{r \in \mathbb{Z} + \frac{1}{2} 
+ \frac{j}{k}}\atop{r>0}}r \ \alpha^{(m), - }(-r)  \alpha^{(m), + }(r)  
 \Biggr) + \frac{j^2d}{2k^2},
\end{multline}
we have that $M_{\sigma_\xi}$ is an ordinary $\sigma_\xi$-twisted $V_{fer} \otimes V_{fer}$-module with graded dimension given by 
\begin{equation}
\mathrm{dim}_q M_{\sigma_\xi}=  q^{-d/24 + j^2d/2k^2} \prod_{{r \in \mathbb{Z} + \frac{1}{2} - \frac{j}{k}}\atop{r>0}} (1+q^{r})^{d} \prod_{{r \in \mathbb{Z} + \frac{1}{2} + \frac{j}{k}}\atop{r>0}}  (1 + q^{r} )^d.
\end{equation}
 
The space $V_{bos} \otimes V_{bos} \otimes M_{\sigma_\xi}$ is a $\sigma_\xi$-twisted $V \otimes V$-module with twisted vertex operators $Y^{\sigma_\xi}(u_1 \otimes v_1 \otimes u_2 \otimes v_2, x) = Y(u_1 \otimes u_2, x) \otimes Y^{\sigma_\xi} (v_1 \otimes v_2, x)$ for $u_1, u_2 \in V_{bos}$ and $v_1, v_2 \in V_{fer}$.   And $V_{bos} \otimes V_{bos} \otimes M_{\sigma_\xi}$ is an ordinary $\sigma_\xi$-twisted $V \otimes V$-module.  In the free case, the $q$-dimension is
\begin{multline}
\mathrm{dim}_q V_{bos} \otimes V_{bos} \otimes M_{\sigma_\xi}\\
=  q^{-d/24 + j^2d/2k^2} \eta(q)^{-2d} \prod_{{r \in \mathbb{Z} + \frac{1}{2} - \frac{j}{k}}\atop{r>0}} (1+q^r)^{d} \prod_{{r \in \mathbb{Z} + \frac{1}{2} + \frac{j}{k}}\atop{r>0}}  (1 + q^r )^d.
\end{multline}
Furthermore, we have
\begin{multline}
J^{\sigma_\xi}(0) = \sum_{m = 1}^d \Biggl( \sum_{{r \in \mathbb{Z} + \frac{1}{2}- \frac{j}{k}}\atop{r>0}} \alpha^{(m), +} (-r) \alpha^{(m), -}(r) \\
-  \sum_{{r \in \mathbb{Z} + \frac{1}{2}+ \frac{j}{k}}\atop{r>0}} \alpha^{(m), -} (-r) \alpha^{(m), +}(r) \Biggr),
\end{multline}
and thus the $p,q$-dimension in the free case is
\begin{multline}\label{p,q-dimension-sigma_xi-twisted}
\mathrm{dim}_{p,q} V_{bos} \otimes V_{bos} \otimes M_{\sigma_\xi}\\
=  q^{-d/24 + j^2d/2k^2} \eta(q)^{-2d} \prod_{{r \in \mathbb{Z} + \frac{1}{2} - \frac{j}{k}}\atop{r>0}} (1+p^{-1}q^r)^{d} \prod_{{r \in \mathbb{Z} + \frac{1}{2} + \frac{j}{k}}\atop{r>0}}  (1 + pq^r )^d.
\end{multline}

\begin{rema}{\em 
Let $L$ be a positive definite lattice of rank $d$, let $V_L$ be the vertex operator superalgebra corresponding to $L$, and let $V_{fer}$ be the fermionic vertex operator superalgebra constructed from $\mathfrak{h} = L \otimes_{\mathbb{Z}} \mathbb{C}$.  Then constructing the $\sigma_\xi$-twisted $V_{fer} \otimes V_{fer}$-module, $M_{\sigma_\xi}$, we have that $V_L \otimes V_L \otimes M_{\sigma_\xi}$ is a $(1 \otimes 1 \otimes \sigma_\xi$-twisted $V_L \otimes V_{fer} \otimes V_L \otimes V_{fer}$-module.  In this case, the $p,q$-dimension in the lattice case is given by (\ref{p,q-dimension-sigma_xi-twisted}) multiplied by $\Theta(L)^2$.  However, the twisted modules for extensions of $\sigma_\xi$ acting nontrivially on the odd component of $V_L$ in the case of an integral (not even) lattice have not yet been constructed.}
\end{rema}

\end{document}